\documentclass[a4paper]{amsart}

\usepackage[margin=3cm]{geometry}

\usepackage{graphicx}
\usepackage{tikz}
\usepackage{tikz-cd}
\usepackage{amssymb}
\usepackage{mathtools}
\usepackage{xcolor}
\usetikzlibrary{decorations.markings}
\usetikzlibrary{arrows}

\usepackage{hyperref}
\usepackage{floatrow}

\theoremstyle{theorem}
\newtheorem{theorem}{Theorem}[section]
\newtheorem{proposition}[theorem]{Proposition}
\newtheorem{lemma}[theorem]{Lemma}
\newtheorem{corollary}[theorem]{Corollary}

\theoremstyle{definition}

\newtheorem{notation}[theorem]{Notation}

\theoremstyle{remark}
\newtheorem{remark}[theorem]{Remark}

\numberwithin{equation}{section}

\newcommand{\LL}{\mathcal{L}}
\newcommand{\K}{\mathcal{K}}
\newcommand{\F}{\mathcal{F}}

\newcommand{\OO}{\mathcal{O}}
\newcommand{\NO}{\mathcal{NO}}
\newcommand{\TT}{\mathcal{T}}
\newcommand{\NT}{\mathcal{NT}}

\newcommand{\Z}{\mathbb{Z}}

\newcommand{\C}{\mathbb{C}}
\newcommand{\N}{\mathbb{N}}

\newcommand{\image}{\operatorname{im}}
\newcommand{\spann}{\operatorname{span}}

\newcommand{\coker}{\operatorname{coker}}

\newcommand{\id}{\operatorname{id}}

\newcommand{\hatimes}{\mathbin{\hat{\otimes}}}

\definecolor{aspurple}{rgb}{0.5,0.1,0.65}

\begin{document}
	
\title[$K$-theory of $C^*$-algebras arising from commuting Hilbert bimodules]{$K$-theory of $C^*$-algebras arising from commuting Hilbert bimodules and invariant ideals}

\author[Astrid an Huef]{Astrid an Huef}
\address[A.\ an Huef]{School of Mathematics and Statistics, Victoria University of Wellington, P.O. Box 600, Wellington 6140, New Zealand.}
\email{astrid.anhuef@vuw.ac.nz}
\author[Abraham Ng]{Abraham C.S.\ Ng}
\address[A.C.S. Ng]{School of Mathematics and Statistics, University of Sydney, NSW 2006, Australia.}
\email{abraham.ng@sydney.edu.au}
\author[Aidan Sims]{Aidan Sims}
\address[A. Sims]{School of Mathematics and Statistics,
University of New South Wales, NSW 2052, Australia.} 
\email{aidan.sims@unsw.edu.au}

	
\date{30 March 2025, revised 9 July 2025}

\thanks{This research was supported  by  Marsden grant 21-VUW-156 of the Royal Society of New  Zealand, and by ARC Discovery--Projects grants DP180100595 and DP220101631. We thank an anonymous referee for their helpful comments.}

	\begin{abstract} We study the $K$-theory of the Cuntz--Nica--Pimsner  $C^*$-algebra of a rank-two product system  that is an extension determined by an invariant ideal of the coefficient algebra. We use a construction of Deaconu and Fletcher that describes the Cuntz--Nica--Pimsner $C^*$-algebra of the product system  in terms of two iterations of Pimsner's original construction of a $C^*$-algebra from a right-Hilbert bimodule. We apply our results to the product system built from  two commuting surjective local homeomorphisms of a totally disconnected space, where the Cuntz--Nica--Pimsner  $C^*$-algebra is  isomorphic to the $C^*$-algebra of the associated rank-two Deaconu--Renault groupoid. We then apply a theorem of Spielberg about stable finiteness of an extension to obtain sufficient conditions for stable finiteness of the  $C^*$-algebra of the Deaconu--Renault groupoid.
\end{abstract}
	
\maketitle

 \section{Introduction}

This paper achieves two main objectives. The first concerns rank-two product systems of right-Hilbert bimodules, in the sense of Fowler, in which the left actions are given by injective homomorphisms into the compacts. We show that the inclusion of a suitably invariant ideal in the coefficient algebra induces a morphism of   exact sequences in $K$-theory for the associated Cuntz--Nica--Pimsner algebras. In particular, when the coefficient algebra $A$ has trivial $K_1$ group, we show that the $K_0$-group of the Cuntz--Nica--Pimsner algebra is explicitly computable in terms of the maps on $K_0(A)$ induced by the $KK$-classes of the generating bimodules in the product system. We also show that an invariant ideal $I$ of $A$ yields a morphism between two such exact sequences. The second objective is to apply this to the $C^*$-algebra of a Deaconu--Renault groupoid with totally disconnected unit space. We obtain a criterion for stable finiteness in terms of stable finiteness of the ideal and quotient corresponding to an open invariant subspace of the unit space, and its closed invariant complement.

We came to this investigation through our previous work on stable finiteness of $C^*$-algebras associated to higher-rank graphs \cite{aHNS}. Still earlier work of the first and third authors \cite{CaHS} showed that for a rank-2 graph whose groupoid admits no nontrivial open invariant sets, stable finiteness of its $C^*$-algebra is characterised by a fairly checkable condition on the adjacency matrices of the graph that comes out of work of Voiculescu \cite{Voiculescu} and Brown \cite{Brown} on AF-embeddable $C^*$-algebras. In \cite{aHNS} we parlayed this into a theorem about stable finiteness of $C^*$-algebras of rank-2 graphs whose groupoids admit exactly one nontrivial open invariant set $U$. Our approach to this was to use Spielberg's theorem \cite{Spielberg1988} about stable finiteness of an extension of a stably finite quotient by a stably finite ideal. In order to invoke Spielberg's theorem, we needed a characterisation of the map in $K$-theory induced by the inclusion of the ideal corresponding to $U$ into the $C^*$-algebra of the $2$-graph. But the $K$-theory computations for $2$-graph $C^*$-algebras available at the time \cite{RobSteg, Evans} proceed via Kasparov's spectral sequence \cite{Kasparov}, and do not give enough information to characterise this map. We navigated this problem in \cite{aHNS} by re-computing the $K$-theory of a $2$-graph $C^*$-algebra by iterated applications of the Pimsner--Voiculescu sequence for crossed products by $\mathbb{Z}$. Though this approach does not provide the explicit description of $K_1(C^*(\Lambda))$ that emerges from the spectral sequence, it has the advantage that naturality of the Pimsner--Voiculescu sequence with respect to $\mathbb{Z}$-equivariant homomorphisms gives the desired description of the map between $K_0$-groups induced by the inclusion $I \hookrightarrow C^*(\Lambda)$ of a gauge-invariant ideal of the $2$-graph $C^*$-algebra.

The computation of $K$-theory in \cite{aHNS} has a second drawback: it involves a great deal of very complicated book-keeping
(take a look, for example, at the diagrams on pages~288, 289~and~293 of \cite{aHNS}). However, an alternative approach presents itself: the $C^*$-algebras of rank-2 Deaconu--Renault groupoids can also be described as Cuntz--Pimsner algebras, as defined by Fowler \cite{Fowler}, of product systems $\mathbf{X}$ over $\mathbb{N}^2$. (In keeping with more recent literature, we call them Cuntz--Nica--Pimsner algebras and denote them $\mathcal{NO}_\mathbf{X}$---new constructions such as Sehnem's covariance algebra \cite{Sehnem} reduce to Fowler's in our setting.) Results of Deaconu \cite{Deaconu}, further developed by Fletcher \cite{FletcherThesis, Fletcher2018}, describe this Cuntz--Nica--Pimsner algebra in terms of two iterations of Pimsner's original construction \cite{Pimsner} of a $C^*$-algebra from a right-Hilbert bimodule, and at each iterate, the coefficient algebra of the bimodule is isomorphic to the $C^*$-algebra of the Deaconu--Renault groupoid for a lower-rank dynamics. So at each stage Pimsner's six-term sequence in $KK$-theory \cite{Pimsner}, which is natural for suitable morphisms of right-Hilbert bimodules, computes the $K$-theory of the Cuntz--Pimsner algebra in terms of that of the coefficient algebra and the map in $K$-theory induced by the class of the right-Hilbert bimodule in $KK$-theory.

In this paper we carry out this program for computing $K$-groups and maps between them via iteration of Pimsner's six-term sequence in greater generality.  In \S\ref{sec:ktheoryfletcher}, we consider  rank-2 product systems $\mathbf{X}$ of right-Hilbert bimodules in which the left action of the coefficient algebra $A$ on each module $X_n$ is implemented by an injective homomorphism into the algebra of generalised compact operators. We identify the \emph{$\mathbf{X}$-invariant} ideals $I$ of $A$ for which $\mathbf{X}I \coloneqq (X_n\cdot I)_{n \in \mathbb{N}^2}$ is a product system of right-Hilbert $I$--$I$-bimodules, and show that there is an inclusion $\phi \colon  \mathcal{NO}_{\mathbf{X}I} \hookrightarrow \mathcal{NO}_\mathbf{X}$ whose range is an ideal of $\mathcal{NO}_\mathbf{X}$. We then work out the details of the naturality of Pimsner's six-term sequence in $KK$-theory to obtain a commuting diagram of short exact sequences describing the $K_0$-groups of $\mathcal{NO}_\mathbf{X}$ and of $\mathcal{NO}_{\mathbf{X}I}$ and the map from the latter to the former induced by $\phi$. Even for $2$-graph $C^*$-algebras $C^*(\Lambda)$, which are Cuntz--Nica--Pimsner algebras of product systems over $\mathbb{N}^2$ of right-Hilbert-$C_0(\Lambda)^0$ bimodules, this significantly simplifies the $K$-theory calculations of \cite{aHNS}.

In \S\ref{sec:ktheoryDR} we consider an action $T$ of $\mathbb{N}^2$ by local homeomorphisms of a totally disconnected locally compact Hausdorff space with an open invariant subset $H \subset \Omega$. We identify the $C^*$-algebra $C^*(G_T)$ of the associated Deaconu--Renault groupoid with the Cuntz--Nica--Pimsner algebra of a product system $X$ over $C_0(\Omega)$, and show that $C_0(H) \subset C_0(\Omega)$ is an $X$-invariant ideal. We can therefore apply our earlier results to obtain a new morphism of short-exact sequences that describes the $K_0$-groups of $C^*(G_T)$ and $C^*(G_{T|_H})$ and the map between them induced by the inclusion $C^*(G_{T|_H}) \hookrightarrow C^*(G_T)$. (We do not explore what our results say about $K_1$-groups, since they do not improve upon the description of \cite[Corollary~7.7]{FarsiKumjianPaskSims}.)

Finally in \S\ref{sec:SFE}, as an application, we use a result of Spielberg \cite[Lemma 1.5]{Spielberg1988} to obtain a sufficient condition for stable finiteness of $C^*(G_T)$ in terms of stable finiteness of $C^*(G_{T|_H})$ and $C^*(G_{T|_{\Omega \setminus H}})$ and of the maps on $C_c(\Omega, \mathbb{N})$ induced by the generators $T_{(1,0)}$ and $T_{(0,1)}$ of the action. This is useful because, for example, if the reductions of $T$ to $H$ and to $\Omega\setminus H$ are minimal (so in particular, if the associated $C^*$-algebras are simple), then stable finiteness of the associated $C^*$-algebras is characterised by the results of \cite{BonickeLi} or \cite{RainoneSims2020}.
	
	\section{Preliminaries}
	
	\subsection{Hilbert bimodules and their associated algebras: Basic definitions and results}\label{ssec:HB}
	
	In the following three subsections, we give a brief summary of Hilbert bimodules and the associated algebras (see \cite{RaeburnWilliams} for example). We also establish a consistent map-naming convention that we highlight in the `notation' environments. Given a $C^*$-algebra $A$, a \emph{right-Hilbert $A$-module $X$} is a right $A$-module equipped with an $A$-valued map $\langle\cdot,\cdot\rangle_A \colon X \times X \to A$, called the \textit{inner product}, which is linear in the second variable and such that for $\xi,\eta \in X$ and $a\in A$,
	\begin{itemize}
		\item $\langle \xi,\xi\rangle_A \ge 0,$ with equality if and only if $\xi =0$;
		\item $\langle \xi, \eta\rangle_A = \langle \eta, \xi\rangle_A^*$;
		\item $\langle \xi,  \eta\cdot a\rangle_A = \langle \xi,\eta\rangle_A a$; and
		\item $X$ is complete with respect to the norm given by $\|\xi\|_A^2 = \|\langle \xi,\xi\rangle_A\|$.
	\end{itemize}
	Often we will drop the subscript $A$ on the inner product and associated norm.
	
A map $T\colon X\to X$ is an \textit{adjointable operator} if there exists a map $T^*\colon X\to X$, called the \textit{adjoint} of $T$, such that $\langle T^*\xi,\eta\rangle = \langle \xi,T\eta\rangle$ for all $\xi, \eta\in X$. If $T$ is adjointable, the adjoint $T^*$ is unique and both $T$ and $T^*$ are automatically bounded $A$-module homomorphisms. The space $\mathcal{L}(X)$ of adjointable operators on $X$ is itself a $C^*$-algebra. Given $\xi,\eta \in X$, the map $\theta_{\xi,\eta} \colon  X\to X$ determined by $\theta_{\xi,\eta}(\zeta) = \xi \cdot \langle\eta,\zeta\rangle_A$ is an adjointable operator with adjoint $\theta_{\xi,\eta}^* = \theta_{\eta,\xi}$. We define $\mathcal{K}(X)$ to be $\overline{\spann}\{\theta_{\xi,\eta} \in \mathcal{L}(X):\xi,\eta \in X\}$ (this is an ideal of $\mathcal{L}(X)$) and call it the algebra of \textit{compact operators} on $X$.
	
Given $C^*$-algebras $B, A$, a \textit{right-Hilbert $B$--$A$-bimodule} $X$ is a right-Hilbert $A$-module together with a left action of $B$ by adjointable operators on $X$ that is implemented by a $*$-homomorphism $\phi \colon  B\to \mathcal{L}(X)$ (that is, $b\cdot x = \phi(b)x$). Thus, for all $\xi,\eta\in X, b\in B$,
	\[\langle b\cdot \xi,\eta\rangle_A = \langle \xi, b^*\cdot \eta\rangle_A.\] Sometimes, for clarity, we will write $X = {}_{B}X_A$. If the homomorphism $\phi$ implementing the left action takes values in $\mathcal{K}(X)$, we say that the Hilbert bimodule $X$ has \textit{compact left action}. If $B=A$, we call $X$ a right-Hilbert $A$-bimodule (also known as a \textit{$C^*$-correspondence over $A$}).
		
	Given two right-Hilbert bimodules ${}_{B}X_A$ and ${}_AY_C$, we can first take the quotient of the algebraic tensor product $X \odot Y$ by the submodule generated by differences of the form $\xi \cdot a \odot \eta - \xi \odot a\cdot \eta$ and then complete it with respect to the inner product that satisfies $\langle \xi \odot \eta,\xi'\odot \eta'\rangle_C = \langle \eta,\langle \xi,\xi'\rangle_A \cdot \eta'\rangle_C$. The result is a new right-Hilbert $B$--$C$-bimodule called the \textit{balanced tensor product} $X\otimes_A Y$ with left action $b\cdot(\xi\otimes\eta) = (b\cdot \xi) \otimes \eta$ and right action $(\xi\otimes \eta)\cdot c = \xi\otimes(\eta\cdot c)$.
	
Given a sequence $(X_n)_n$ of right-Hilbert $A$-modules, we define
\[
\bigoplus_{n=0}^\infty X_n := \{(x_n)_n \in \prod_n X_n : \sum_n \langle x_n, x_n\rangle_A\text{ converges in $A$}\};
\]
this is again a right-Hilbert $A$-module with $(x_n)_n \cdot a = (x_n \cdot a)_n$ and $\langle (x_n)_n, (y_n)_n\rangle_A = \sum_n \langle x_n, y_n\rangle_A$. If each $X_n$ is a Hilbert $A$--$A$-bimodule, then so is $\bigoplus_n X_n$ with left action $a \cdot (x_n)_n = (a \cdot x_n)_n$.

Every $C^*$-algebra $A$ can be regarded as a right-Hilbert $A$-bimodule over itself with left and right actions given by honest multiplication and inner product defined by $\langle a,b\rangle = a^*b$. This is denoted by ${}_{A}A_A$. Let $X$ be a right-Hilbert $A$-bimodule. For $n\ge 1$, we write
	\[X^{\otimes n} = \underbrace{X\otimes_AX\otimes_A \cdots \otimes_A X}_{n \text{ terms}},\] and after setting $X^{\otimes 0}$ to be $A$, the space
	\[\mathcal{F}_X = \bigoplus_{n=0}^\infty X^{\otimes n}\] is called the \textit{Fock space} of $X$ and is itself a right-Hilbert $A$-bimodule.
	
	Given two right-Hilbert bimodules ${}_{A} X_A$ and ${}_B Y_B$, a pair of maps $(\lambda\colon  A\to B, \varphi\colon X\to Y)$ is a \textit{right-Hilbert bimodule morphism} if $\lambda$ is a $*$-homomorphism and for all $a\in A, \xi, \eta \in X$,
	\[\varphi(a\cdot \xi) = \lambda(a)\cdot\varphi(\xi), \quad \varphi(\xi\cdot a) = \varphi(\xi)\cdot\lambda(a), \quad \text{ and }\quad \lambda (\langle \xi,\eta\rangle) = \langle \varphi(\xi),\varphi(\eta)\rangle.\]
	
	Given a right-Hilbert bimodule ${}_A X_A$, a \textit{Toeplitz representation} of $X$ in a $C^*$-algebra $B$ is a right-Hilbert bimodule morphism $(\pi\colon  X\to B, \psi \colon  X\to {}_B B_B)$. There exists a universal Toeplitz representation of $X$ in a $C^*$-algebra which we denote $\mathcal{T}_X$ and call the \textit{Toeplitz algebra} of $X$. We will denote the representation in the Toeplitz algebra by
	\[(i_1 \colon A \to \mathcal{T}_X, i_2\colon  X\to \mathcal{T}_X).\]
We call $A$ the \textit{coefficient algebra} of $\mathcal{T}_X$. The map $i_1$ is injective and $\TT_X$ is generated by the images of $i_1$ and $i_2$.
	
	\begin{notation}\label{notation:i_1andi_2}
		Given a right-Hilbert $A$-bimodule $X$, the symbol $i_1$ always denotes the canonical inclusion of $A$ into $\TT_X$ and $i_2$ always denotes the canonical map taking $X$ into $\TT_X$.
	\end{notation}
	
	By \cite[Proposition~3.3]{Pimsner}, we can identify $\mathcal{T}_X$ with a subalgebra of $\mathcal{L}(\mathcal{F}_X)$ such that given $a\in A$ and $\xi\in X$,
		\[i_1(a) \colon (\xi^n)_n^{\infty} \in \mathcal{F}_X\mapsto (a\cdot \xi^n)_{n=0}^\infty \in \mathcal{F}_X,\] and
		\[i_2(\xi) \colon (\xi^n)_{n=0}^{\infty} \in \mathcal{F}_X\mapsto (0,\xi\cdot \xi^0,\xi\otimes \xi^1,\xi\otimes\xi^2,\xi\otimes\xi^3,\dots) \in \mathcal{F}_X.\]
	Furthermore, we have that $\mathcal{K}(\mathcal{F}_X) \lhd \TT_X \subset \mathcal{L}(\mathcal{F}_X)$.
	
	Given any Toeplitz representation $(\pi, \psi)$ of a right-Hilbert bimodule ${}_A X_A$ in a $C^*$-algebra $B$, there exists a $*$-homomorphism $(\pi,\psi)^{(1)}\colon \mathcal{K}(X)\to B$ such that for all $\xi,\eta\in X$,
	\[(\pi,\psi)^{(1)}(\theta_{\xi,\eta}) = \psi(\xi)\psi(y)^*.\]
As usual, for an ideal $I \lhd A$ we write $I^\perp := \{a\in A : aI = \{0\}\}$.
We say that a Toeplitz representation $(\pi, \psi)$ is \textit{Cuntz--Pimsner covariant} if $(\pi,\psi)^{(1)}(\phi(a)) = \pi(a)$ for all $a\in \phi^{-1}(\mathcal{K}(X)\cap (\ker\phi)^\perp)$. There exists a universal Cuntz--Pimsner covariant Toeplitz representation of $X$ in a $C^*$-algebra which we denote $\mathcal{O}_X$ and call the \textit{Cuntz--Pimsner algebra} of $X$. We will denote the representation in the Cuntz--Pimsner algebra by
	\[(j_1 \colon A \to \mathcal{O}_X, j_2\colon  X\to \mathcal{O}_X).\]
As with the Toeplitz representation, $j_1$ is always injective and $\OO_X$ is generated by the images of $j_1$ and $j_2$. If the left action homomorphism $\phi$ is injective, then $\OO_X$ is the quotient of $\mathcal{T}_X$ by the ideal of $\TT_X$ generated by $\{\pi(a) - (i_1,i_2)^{(1)}(\phi(a)) : \phi(a) \in \mathcal{K}(X)\}$.	In general, $\OO_X$ is the quotient of $\mathcal{T}_X$ by a sub-ideal of $\mathcal{K}(\mathcal{F}_X)$ and the quotient map $q\colon \TT_X\to \OO_X$ satisfies $j_1 = q\circ i_1$ and $j_2 = q\circ i_2$.
	
	\begin{notation}
		Given a right-Hilbert $A$-bimodule $X$, the homomorphism $j_1$ will always denote the canonical inclusion of $A$ into $\OO_X$, $j_2$ will always denote the canonical map taking $X$ to $\OO_X$, and $q$ will always denote the quotient map $\TT_X\to \OO_X$.
	\end{notation}
	
	Any right-Hilbert $A$-bimodule morphism $(\lambda \colon  A\to B, \varphi\colon  {}_A X_A \to {}_B Y_B)$ induces three important homomorphisms:
	\begin{enumerate}
	\item  $\Phi\colon  \F_X\to \F_Y$ such that $\xi_1\otimes\cdots\otimes \xi_n \in X^{\otimes n} \mapsto \varphi(\xi_1)\otimes\cdots \otimes\varphi(\xi_n) \in Y^{\otimes n}$ (this induces an additional homomorphism $\Phi^{(1)}\colon  \K(\F_X)\to\K(\F_Y)$ such that $\theta_{\bar{\xi},\bar{\eta}}\mapsto \theta_{\Phi(\bar{\xi}),\Phi(\bar{\eta})}$ for $\bar{\xi},\bar{\eta} \in \F_X$).
	\item $\TT_X \to \TT_Y$ such that $i_1(a) \mapsto i_1(\lambda(a))$ and $i_2(\xi) \mapsto i_2(\varphi(\xi))$.
    \item $\varphi^{(1)} : \K(X) \to \K(Y)$ such that $\theta_{\xi, \eta} \mapsto \theta_{\varphi(\xi), \varphi(\eta)}$ (this is the restriction of $\Phi^{(1)}$ from~(1) to $\K(X) \subset \K(\F_X)$).
    \end{enumerate}
    If $(\lambda, \varphi)$ is \emph{covariant} (also known as \emph{coisometric} \cite[Definition~1.3]{Brenken}) in the sense that $\varphi^{(1)}(\phi(a)) = \varphi(a)$ whenever $\phi(a) \in \K(X)$, and $\lambda$ is injective, then by \cite[Corollary~1.5]{Brenken} there is a fourth homomorphism
    \begin{enumerate}\setcounter{enumi}{3}
    \item $\OO_X \to \OO_Y$ such that $j_1(a) \mapsto j_1(\lambda(a))$ and $j_2(\xi) \mapsto j_2(\varphi(\xi))$ for all $a\in A, \xi\in X$.
	\end{enumerate}
	
	Given a right-Hilbert bimodule ${}_A X_{A}$, an ideal $I \lhd A$ is \textit{$X$-invariant} if $I X \subset X I$ (here and elsewhere, $IX$ means $I\cdot X$ and $XI$ means $X\cdot I$). For any $X$-invariant ideal $I$, ${}_I X I_I$ is a right-Hilbert bimodule and the pair $(I\hookrightarrow A, XI \hookrightarrow X)$ is a coisometric right-Hilbert bimodule morphism. The induced homomorphism $\OO_{XI} \to \OO_X$ is injective (see \cite[Corollary 3.9]{FMR}).

	\subsection{Hilbert bimodules, the associated algebras and \texorpdfstring{$K$}{K}-theory}\label{ssec:HBktheory}
	
	Let $X$ be a right-Hilbert $A$-bimodule. Then $i_1 \colon  A\to \TT_X$ induces a $KK$-equivalence so that $(i_1)_* \colon  K_0(A) \to K_0(\TT_X)$ is an isomorphism \cite[\S4]{Pimsner}.
	
	Let $P \in \mathcal{K}(\mathcal{F}_X)$ be the infinite matrix such that $P_{11} = 1$ and $P_{ij}=0$ whenever $i\neq 0$ or $j\neq 0$. Let $k\colon A\to \mathcal{K}(\mathcal{F}_X)$ be the map $a\mapsto aP$. Then for all $a\in A$,
	\[k(a)\colon  (\xi^n)_{n=0}^\infty \mapsto (a\xi^0,0,0,\dots).\] It is not hard to see that $\mathcal{K}(\mathcal{F}_X)P\mathcal{K}(\mathcal{F}_X) = \mathcal{K}(\mathcal{F}_X)$ and that $P\mathcal{K}(\mathcal{F}_X)P = k(A) \cong A$. Thus $k_* \colon  K_0(A) \to K_0(\mathcal{K}(\mathcal{F}_X))$ is an isomorphism.

	\begin{notation}
		Given a right-Hilbert bimodule $X$ over a $C^*$-algebra $A$, the letter $k$ always denotes the map $A\to \mathcal{K}(\mathcal{F}_X), \;a \mapsto aP$, and $\ell$ always denotes the inclusion $\mathcal{K}(\mathcal{F}_X)\hookrightarrow \TT_X$.
	\end{notation}

	There is a 6-term exact $K$-theory sequence called the Pimsner sequence associated to any right-Hilbert bimodule $X$ over a $C^*$-algebra $A$ with left action by compact operators. 	Using the maps $i_1,j_1,k,\ell$ of our map-naming convention, the Pimsner sequence can then be written as the outer sequence in the following diagram:
	\[
	\begin{tikzcd}
		K_0(A) \arrow[d, "{\cong ,k_*}"] \arrow[r, "{1-[X]}"]                          & K_0(A) \arrow[d, "{\cong,(i_1)_*}"] \arrow[rd, "(j_1)_*", dashed, bend left]  &                                                                                \\
		K_0(\K(\F_X)) \arrow[r, "\ell_*"]                                              & K_0(\TT_X) \arrow[r, "q_*"]                                                   & K_0(\OO_X) \arrow[d, "\partial"] \arrow[dd, "\partial", dashed, bend left=70] \\
		K_1(\OO_X) \arrow[u, "\partial"] \arrow[uu, "\partial", dashed, bend left=70] & K_1(\TT_X) \arrow[l, "q_*"]                                                   & K_1(\K(\F_X)) \arrow[l, "\ell_*"]                                              \\
		& K_1(A) \arrow[u, "{\cong, (i_1)_*}"] \arrow[lu, "(j_1)_*", dashed, bend left] & K_1(A) \arrow[u, "{\cong,k_*}"] \arrow[l, "{1-[X]}"]
	\end{tikzcd}.\]
	We will work with both the inner 6-term sequence and the outer 6-term sequence, sometimes in $K$-theory and sometimes in $KK$-theory, depending on the situation.  The instances of $[X]$ in the diagram above denote the  homomorphisms $[X]\colon K_*(A)\to K_*(A)$  obtained from the composition
	\begin{equation}\label{eq:KKnotation}
 \begin{tikzcd}
	K_*(A) \arrow[r, "{\Theta_*}"]
	& {KK_*(\C,A)} \arrow[r, "\cdot\hatimes X"] & {KK_*(\C,A)} \arrow[r, "{\Theta_*^{-1}}"] & K_*(A)
\end{tikzcd}
\end{equation}
using the isomorphisms $\Theta\colon K_*(A)\to KK_*(\C,A)$  and  the Kasparov products $\cdot\hatimes X$.  In our treatment below, $A$ is approximately finite-dimensional, and hence $K_1(A)=0$. In preparation for our application to the $K$-theory of a rank-2 Deaconu--Renault groupoid,  in \S\ref{sec:KK} we discuss $[X]$ in detail when $X$ is the right-Hilbert module associated to a surjective local homeomorphism between totally disconnected spaces.

	\subsection{Product systems and the associated algebras}\label{ssec:prodsys}
	
	Apart from a brief appearance in \S\ref{sec:ktheoryDR}, the main part of this paper can be read without reference to product systems, which were first introduced in \cite{Fowler}. This preliminary subsection is mostly taken from \cite[\S3.1]{FletcherThesis} and included for those who want to understand the technical proofs in \S\ref{ssec:appendix1}  concerning the Deaconu--Fletcher constructions.
	
	Let $A$ be a $C^*$-algebra and let $P$ a countable semigroup with identity $e$. A \textit{product system over $P$ with coefficient algebra $A$} (as defined in \cite[Definition 3.1.1]{FletcherThesis}) is a semigroup $\mathbf{X} = \bigsqcup_{p\in P} \mathbf{X}_p$ such that for $p,q \in P$,
	\begin{itemize}
		\item $\mathbf{X}_p\subset \mathbf{X}$ is a right-Hilbert $A$-bimodule;
		\item $\mathbf{X}_0$ is equal to  ${}_A A_A$;
		\item For $p,q \in P$ with $p \not= 0$, there exists a right-Hilbert bimodule isomorphism $M_{p,q} \colon  \mathbf{X}_p\otimes_A \mathbf{X}_q \to \mathbf{X}_{pq}$
            satisfying $M_{p,q}(\xi\otimes \eta) = \xi\eta$ for all $\xi \in \mathbf{X}_p$ and $y\in \mathbf{X}_q$;
		\item For $a\in \mathbf{X}_0 = A$ and $\xi \in \mathbf{X}_p$, $a\xi = a\cdot \xi$ and $\xi a = \xi \cdot a$.
	\end{itemize}
	
	Product systems generalise right-Hilbert bimodules: if $X$ is a  right-Hilbert $A$-bimodule, then $\mathbf{X} = \bigsqcup_{n\in \N} X^{\otimes n}$ is a product system.
	
	In our paper, we only use  $P=\N^2$, and so we restrict the background material to this case. Furthermore, we will always assume that the left actions are by compacts (implying that the product system is \textit{compactly aligned}, though we will not delve further into this term). A \textit{representation} of $\mathbf{X}$ in a $C^*$-algebra $B$ is a map $\psi\colon  \mathbf{X} \to B$ such that for $p,q\in \N^2$,
	\begin{itemize}
		\item $\psi_p \coloneqq\psi|_{\mathbf{X}_p}$ is a linear map;
		\item $\psi_0$ is a $*$-homomorphism;
		\item $\psi_p(\xi)\psi_q(\eta) = \psi_{pq}(\xi\eta)$ for all $\xi\in \mathbf{X}_p, \eta\in \mathbf{X}_q$; and
		\item $\psi_p(\xi)^*\psi_p(\eta) = \psi_0(\langle\xi,\eta\rangle_A)$ for all $\xi,\eta\in\mathbf{X}_p$.
	\end{itemize}
	Note that for all $p\in \N^2$, $(\psi_0,\psi_p)$ is a Toeplitz representation of $\mathbf{X}_p$ in $B$.
	
	For $p=(p_1,p_2),q=(q_1,q_2) \in \N^2$, let $p\vee q = (\max\{p_1,q_1\},\max\{p_2,q_2\})\in \N^2$. For $S\in \K(\mathbf{X}_p), T\in \K(\mathbf{X}_q)$, let
	\[
        \iota_{p}^{p\vee q}(S) = M_{p,(p\vee q)-p}\circ (S\otimes \operatorname{id}_{\mathbf{X}_{p\vee q - p}})\circ M_{p,(p\vee q)-p}^{-1} \in \K(\mathbf{X}_{p\vee q}),
    \] 
    and
	\[
        \iota_{q}^{p\vee q}(T) = M_{q,(p\vee q)-q}\circ (T\otimes \operatorname{id}_{\mathbf{X}_{p\vee q - q}})\circ M_{q,(p\vee q)-q}^{-1} \in \K(\mathbf{X}_{p\vee q}).
    \]
	If $\psi$ satisfies the additional condition that for all $p,q\in\N^2$,
	\begin{itemize}
		\item $(\psi_0,\psi_p)^{(1)}(S)(\psi_0,\psi_q)^{(1)}(T) = (\psi_0,\psi_{p\vee q})^{(1)}(\iota_{p}^{p\vee q}(S)\iota_q^{p\vee q}(T))$ for all $S\in \K(\mathbf{X}_p), T\in \K(\mathbf{X}_q)$,
	\end{itemize}
	then the representation is said to be \textit{Nica covariant}. There exists a univeral Nica covariant respresentation of $\mathbf{X}$ in a $C^*$-algebra which we denote $\mathcal{NT}_\mathbf{X}$ and call the \textit{Nica--Toeplitz algebra} of $\mathbf{X}$. We will always write the representation in the Nica--Toeplitz algebra as
	\[i_\mathbf{X} \colon  \mathbf{X}\to \mathcal{NT}_\mathbf{X},\] and we write $i_{\mathbf{X}_p} = i_{\mathbf{X}}|_{\mathbf{X}_p}$. The algebra $\mathcal{NT}_\mathbf{X}$ is generated by the image of $i_\mathbf{X}$.
	
	There is a further technical condition that we will not delve into that makes a Nica covariant representation what we call \textit{Cuntz--Pimsner covariant}. Suffice it to say that if the left actions are all faithful, something we will always assume, then there exists a universal Cuntz--Pimsner covariant representation of $\mathbf{X}$ in a $C^*$-algebra which we denote $\mathcal{NO}_{\mathbf{X}}$ and call the \textit{Cuntz--Nica--Pimsner algebra} of $\mathbf{X}$. We will always write the representation in the Cuntz--Nica--Pimsner algebra as
	\[j_{\mathbf{X}}\colon  \mathbf{X} \to \mathcal{NO}_{\mathbf{X}},\] and we write $j_{\mathbf{X}_p} = j_{\mathbf{X}}|_{\mathbf{X}_p}$
	The algebra $\NO_{\mathbf{X}}$ is generated by $j_\mathbf{X}(X)$ and is a quotient of $\NT_{\mathbf{X}}$. We denote the quotient map by $q_{\mathbf{X}}\colon \NT_\mathbf{X} \to \NO_{\mathbf{X}}$. We have $j_\mathbf{X} = q_\mathbf{X} \circ i_\mathbf{X}$.
	
\section{The Deaconu--Fletcher construction for commuting Hilbert bimodules}\label{sec:fletcher}

	Let $A$ be a separable $C^*$-algebra, and
	let $X_1$ and $X_2$ be countably generated right-Hilbert $A$-bimodules   with faithful compact left actions. Suppose that there is an isomorphism $X_1 \otimes_A X_2 \to X_2\otimes_A X_1$; we then say that $X_1$ and $X_2$ are \textit{commuting} Hilbert bimodules.  Deaconu and Fletcher produced several constructions which we now describe (see \cite[Propositions~4.2.11 and~4.3.1]{FletcherThesis}).
	
    If we restrict the left action of ${_{\mathcal{T}_{X_1}}\!(\mathcal{T}_{X_1})_{\mathcal{T}_{X_1}}}$ to $A$ (strictly speaking, it is to $i_1(A)$), then we get a right-Hilbert bimodule ${}_{A}(\mathcal{T}_{X_1})_{\mathcal{T}_{X_1}}$. Thus we can take the balanced tensor product $X_2 \otimes_A \mathcal{T}_{X_1}$ which is a right-Hilbert  $A$--$\TT_{X_1}$-bimodule. There is a left action by $\TT_{X_1}$ on $X_2 \otimes_A \mathcal{T}_{X_1}$ defined as follows: the homomorphism of $A \to \mathcal{L}(X_2 \otimes_A \mathcal{T}_{X_1})$ determined by the left action of $A$ on $X_2$ and the linear map $X_1 \to \mathcal{L}(X_2 \otimes_A \mathcal{T}_{X_1})$ given by the multiplication isomorphism $X_1 \otimes X_2 \to X_2 \otimes_A X_1$ constitute a Toeplitz representation of $X_1$ in $\mathcal{L}(X_2 \otimes_A \mathcal{T}_{X_1})$, and therefore extend to a homomorphism $\mathcal{T}_{X_1} \to \mathcal{L}(X_2 \otimes_A \mathcal{T}_{X_1})$ by the universal property of $\mathcal{T}_{X_1}$ \cite{Deaconu, FletcherThesis}. This gives a right-Hilbert $\TT_{X_1}$-bimodule. The same construction with $X_1$ and $X_2$ switched yields a right-Hilbert $\TT_{X_2}$-bimodule $X_1\otimes_A \TT_{X_2}$. 	
	Similarly, there is a left action by $\OO_{X_1}$ on the right-Hilbert $A$--$\OO_{X_1}$-bimodule $X_2 \otimes_A \OO_{X_1}$ that extends the original left action of $A$. This gives a  right-Hilbert $\OO_{X_1}$-bimodule $X_2 \otimes_A \OO_{X_1}$. The same construction with $X_1$ and $X_2$ switched yields a right-Hilbert $\OO_{X_2}$-bimodule $X_1\otimes_A \OO_{X_2}$.
	
	For $m,n\in\N$, let $\mathbf{X}_{(m,n)} = X_1^{\otimes m}\otimes_A X_2^{\otimes n}$. Since $X_1\otimes_A X_2$ is isomorphic to $X_2\otimes_A X_1$, $\mathbf{X} = \bigsqcup_{(m,n)\in\N^2} \mathbf{X}_{(m,n)}$ is a product system with associated Nica--Toeplitz and Cuntz--Nica--Toeplitz algebras (see \S\ref{ssec:prodsys}).
	
	In the next theorem we summarise theorems from the Deaconu--Fletcher construction \cite[Theorems~3.3.17, 3.4.21, 4.2.12, 4.3.2, 4.3.3]{FletcherThesis}, as they pertain to our situation.
		
	\begin{theorem}[Deaconu--Fletcher]\label{thm:flecther}
		Let $X_1,X_2$ be as above. Then there are isomorphisms (see Remark~\ref{rmk:expfletcheriso} below)
		\begin{enumerate}
			\item[(1)] $\TT_{X_2 \otimes_A \TT_{X_1}} \cong \TT_{X_1\otimes_A \TT_{X_2}} \cong \NT_\mathbf{X}$;
			\item[(2)] $\OO_{X_2 \otimes_A \OO_{X_1}} \cong \OO_{X_1\otimes_A \OO_{X_2}} \cong \NO_\mathbf{\mathbf{X}}$;
			\item[(3)] $\OO_{X_2 \otimes_A \TT_{X_1}} \cong \TT_{X_1\otimes_A \OO_{X_2}}$ and $\OO_{X_1 \otimes_A \TT_{X_2}} \cong \TT_{X_2\otimes_A \OO_{X_1}}$.
		\end{enumerate}
	\end{theorem}

	\begin{remark}\label{rmk:expfletcheriso}
		In particular, in the notation from \S\ref{ssec:HB}, the isomorphism $\OO_{X_2\otimes_A\OO_{X_1}}\to\NO_{\mathbf{X}}$ does the following:  for $a\in A$, $\xi \in X_1$ and $\eta\in X_2$,
		\begin{gather*}
		j_1(j_1(a))\mapsto j_{\mathbf{X}_{(0,0)}}(a), \quad j_1(j_2(\xi))\mapsto j_{\mathbf{X}_{(1,0)}}(\xi),\\
		 j_2(\eta\otimes j_1(a))\mapsto j_{\mathbf{X}_{(0,1)}}(\eta\cdot a),\quad  j_2(\eta\otimes j_2(\xi)) \mapsto j_{\mathbf{X}_{(1,1)}}(\eta\xi).\end{gather*}
    The formulas for the isomorphisms in (1)~and~(3) are similar, but with the universal maps $i_*$ into Toeplitz algebras replacing the universal maps $j_*$ into Cuntz--Pimsner algebras as appropriate; the details for~(3) are spelled out in Lemma~\ref{lem:fletcherisomorphism}, and for~(1), all instances of $j$ are replaced by $i$ in the formulas above.
	\end{remark}
	
		We now state four lemmas regarding the Deaconu--Fletcher construction that will be used to compute $K$-theory in \S\ref{sec:ktheoryfletcher}. The proofs are straightforward but very technical and thus are provided in the appendix.
	
	\begin{lemma}\label{lem:fcbimodulemorphism}
		Let $k\colon  A\to \K(\F_{X_2}), \ell \colon  \K(\F_{X_2})\to \TT_{X_2}$ be the homomorphisms associated to the right-Hilbert $A$-bimodule $X_2$. Then there is a map $\psi\colon X_1\to X_1\otimes_A \TT_{X_2}$ such that  $\psi(\xi\cdot a) = \xi\otimes \ell(k(a))$ for $a\in A$ and $\xi\in X$. The pair $(\ell\circ k,\psi) \colon  (A,X_1)\to (\TT_{X_2}, X_1\otimes_A \TT_{X_2})$ is a covariant right-Hilbert bimodule morphism.
	\end{lemma}
	
	Lemma~\ref{lem:fcbimodulemorphism} is proved in Appendix~\ref{sec:appendix} on page~\pageref{proof_lem:fcbimodulemorphism}.
	
	\begin{lemma}\label{lem:incbimodulemorphism}
		Let $i_1 \colon A\to \TT_{X_2}$ be the homomorphism associated to the right-Hilbert $A$-bimodule $X_2$.
		Then there is a map  $\varphi\colon X_1\to X_1\otimes_A \TT_{X_2}$ such that $\varphi(\xi\cdot a) = \xi\otimes i_1(a)$ for $a\in A$ and $\xi\in X_1$. The pair $(i_1,\varphi) \colon (A,X_1)\to (\TT_{X_2}, X_1\otimes_A \TT_{X_2})$ is a covariant right-Hilbert bimodule morphism.
	\end{lemma}
	
	We prove Lemma~\ref{lem:incbimodulemorphism} in Appendix~\ref{sec:appendix}  on page~\pageref{proof_lem:incbimodulemorphism}.
	
	\begin{lemma}\label{lem:fletcherisomorphism}
		Let $J\colon   \OO_{X_1\otimes_A \TT_{X_2}} \to \TT_{X_2\otimes_A \OO_{X_1}}$ be the isomorphism from Theorem~\ref{thm:flecther}(3). Adopt the notation of \S\ref{ssec:HB} so that $(i_1,i_2)$ denotes the Toeplitz representations of both $(A,X_2)$ and $(\OO_{X_1}, X_2\otimes_A \OO_{X_1})$ and $(j_1,j_2)$ denotes the covariant Cuntz--Pimsner representations of both $(\TT_{X_2},X_1\otimes_A\TT_{X_2})$ and $(A,X_1)$.  Fix $a \in A$, $\xi\in X_1$ and  $\eta \in X_2$, and write $\eta = \eta'\cdot \langle\eta',\eta'\rangle \in X_2$ for $\eta' \in X_2$ using \cite[Proposition~2.33]{RaeburnWilliams}. Then
		\begin{gather*}
            J(j_2(\xi\otimes i_1(a))) = i_1(j_2(\xi\cdot a)), \qquad  J(j_1(i_1(a))) = i_1(j_1(a)),\qquad\text{ and }\\
             J(j_1(i_2(\eta))) = i_2(\eta' \otimes j_1(\langle \eta',\eta'\rangle)).
        \end{gather*}
	\end{lemma}
	We prove Lemma~\ref{lem:fletcherisomorphism} in Appendix~\ref{sec:appendix}  on page~\pageref{proof_lem:fletcherisomorphism}.

	Recall that for a right-Hilbert $A$-bimodule $X$, an ideal $I\lhd A$ is \emph{$X$-invariant} if $IX\subset XI$. Let $I$ be an $X_1$- and $X_2$-invariant ideal of $A$. Then for $i=1,2$, $X_i I$ is a  right-Hilbert  $I$-bimodule and $X_i I \subset X_i$. By \cite[Corollary 3.9]{FMR}, the canonical map $\kappa_i \colon \OO_{X_i I} \to \OO_{X_i}$ is injective. Thus we have the following lemma.
	
	\begin{lemma}\label{lem:fletchermorphisms}
		Let $I$ be an $X_1$- and $X_2$-invariant ideal of $A$. Then
		\begin{enumerate}
			\item
			The pair $(I\hookrightarrow A, X_1 I\hookrightarrow X_1)$ is a covariant right-Hilbert bimodule morphism; and
			\item
			The pair $(\kappa_1 \colon  \OO_{X_1 I}\to \OO_{X_1}, 1\otimes \kappa_1 \colon  X_2I\otimes_I \OO_{X_1 I} \to X_2 \otimes_A \OO_{X_1})$ is a right-Hilbert bimodule morphism. In particular, this induces a homomorphism $\phi\colon  \OO_{X_2I\otimes_I \OO_{X_1 I}}\to \OO_{X_2 \otimes_A \OO_{X_1}}$.
		\end{enumerate}
	\end{lemma}
	We prove Lemma~\ref{lem:fletchermorphisms} in Appendix~\ref{sec:appendix} on page~\pageref{proof_lem:fletchermorphisms}.

\section{\texorpdfstring{$K$}{K}-theory for Cuntz--Pimsner algebras of commuting Hilbert bimodules}\label{sec:ktheoryfletcher}

In this section we prove a  generalisation of \cite[Theorem 4.20]{aHNS} for general commuting right-Hilbert  $A$-bimodules $X_1$ and $X_2$, both with compact left actions and   $K_1(A)=0$, and an ideal $I$ of $A$ which is $X_1$- and $X_2$-invariant.  While the conclusion of Theorem~\ref{thm:naturalses} looks very much like  \cite[Theorem 4.20]{aHNS}, even our application of the theorem to rank-$2$ Deaconu--Renault  groupoids yields new results (see Theorem~\ref{thm:DRkgroups2} below). Our proof is based on Theorem~\ref{thm:diagram1}. Consequently our method of proof is very different from the one for   \cite[Theorem 4.20]{aHNS}.  Notice that the Cuntz--Nica--Pimsner $C^*$-algebra $\mathcal{NO}_\mathbf{X}$ of the product system $\mathbf{X}$ induced by $X_1$ and $X_2$, which we discussed in the introduction,  is  isomorphic by Theorem~\ref{thm:flecther} to the Cuntz--Pimsner algebra  of the bimodule $Y_1\coloneqq X_2 \otimes_A \OO_{X_1}$ appearing in the theorem.

	\begin{theorem}\label{thm:diagram1} Let $A$ be a separable $C^*$-algebra such that $K_1(A)=0$.
		Let $X_1$ and $X_2$ be countably generated  right-Hilbert $A$-bimodules  with faithful compact left actions  of $A$ and such that $X_1\otimes_A X_2 \cong X_2\otimes_A X_1$. 		Let $Y_1$ denote the  right-Hilbert $\OO_{X_1}$-bimodule $X_2 \otimes_A \OO_{X_1}$ from the Deaconu--Fletcher construction. Then  the following diagram commutes:
\begin{equation}\label{eq:monster diagram}
\begin{tikzpicture}[baseline=(a)]
			\node (a) at (0,0){
				\begin{tikzcd}
					&[2em] K_1(\OO_{Y_1}) \arrow[ddl, "\partial", bend right=25, out=320]  \arrow[ddrrr, bend left=19, leftarrow, swap, pos=0.3, "(j_1)_*"]                                           &                                             & K_0(A) \arrow[r, pos=0.55, "{1-[X_2]}"]                               &[2em] K_0(A)                            \\
					0                                                           & 0                                             &                                                        & K_0(A) \arrow[r, crossing over, pos=0.55, "{1-[X_2]}"] \arrow[u, crossing over, pos=0.4, "{1-[X_1]}"]        & K_0(A) \arrow[u, pos=0.4, "{1-[X_1]}"]       \\
					K_0(\OO_{X_1}) \arrow[u, shift left] \arrow[r, "{1-[Y_1]}"] & K_0(\OO_{X_1}) \arrow[u] \arrow[r, "(j_1)_*"] & K_0(\OO_{Y_1}) \arrow[r, "\partial"] & K_1(\OO_{X_1}) \arrow[u, "\partial"] \arrow[r, pos=0.55, "{1-[Y_1]}"] & K_1(\OO_{X_1}) \arrow[u, "\partial"]\\
					K_0(A) \arrow[u, "(j_1)_*"] \arrow[r, "{1-[X_2]}"]          & K_0(A) \arrow[u, "(j_1)_*"]                   &                                                        & 0 \arrow[u]                                                 & 0 \arrow[u]                         \\
					K_0(A) \arrow[r, "{1-[X_2]}"] \arrow[u, "{1-[X_1]}"]        & K_0(A). \arrow[u, "{1-[X_1]}"]                 &                                                        &                                                             &
				\end{tikzcd}
			};
		\end{tikzpicture}\end{equation}
Furthermore, if $I$ is an $X_1$- and $X_2$-invariant ideal of $A$,  then $X_1I$ and $X_2I$ commute, and with  $Y_1^I \coloneqq (X_2 I)\otimes_I \OO_{X_1 I}$ the diagram is natural  with respect to the right-Hilbert bimodule morphisms from Lemma \ref{lem:fletchermorphisms}.
	\end{theorem}
	
	To prove that the diagram~\eqref{eq:monster diagram} commutes, and in particular that the square on the right side involving the two maps $\partial$ commutes, we will use the following proposition and Lemma~\ref{lem:compositionfrom Pimsner} below.
	
	\begin{proposition}\label{prop:square3} Let $(\ell\circ k,\psi) \colon  (A,X_1)\to (\TT_{X_2}, X_1\otimes_A \TT_{X_2})$ be the morphism of right- Hilbert bimodules from Lemma~\ref{lem:fcbimodulemorphism}, and let \[\eta\colon \K(\F_{X_1})\to \K(\F_{X_1\otimes \TT_{X_2}})\text{\ and\ }\theta \colon  \OO_{X_1}\to \OO_{X_1\otimes \TT_{X_2}}\] be the induced homomorphisms. Let $(i_1,\varphi) \colon  (A,X_1)\to (\TT_{X_2}, X_1\otimes_A \TT_{X_2})$ be the morphism of  right-Hilbert bimodules from Lemma~\ref{lem:incbimodulemorphism},
and let \[\rho\colon \K(\F_{X_1})\to \K(\F_{X_1\otimes \TT_{X_2}}) \text{\ and\ } \sigma\colon  \OO_{X_1}\to \OO_{X_1\otimes \TT_{X_2}}\] be the induced homomorphisms. Finally, let $J\colon \OO_{X_1\otimes \TT_{X_2}}\to \TT_{X_2\otimes \OO_{X_1}}$ be the isomorphism of Lemma \ref{lem:fletcherisomorphism}. In the diagram
\[
\begin{tikzpicture}[scale=0.9,baseline=(a), every node/.style={scale=0.9}]
\node (a) at (0,0){
\begin{tikzcd}[column sep = small]
   {K_0(\K(\F_{X_1}))} \arrow[rrrrrdd, "\eta_*"] 
 & 
 & K_0(A) \arrow[ll, swap, "{\cong,k_*}"] \arrow[rr, "{\cong,k_*}"] 
 & 
 & {K_0(\K(\F_{X_2}))} \arrow[r, "\ell_*"] \arrow[rdd, pos=0.4, "(a)", phantom] 
 & {K_0(\TT_{X_2})} \arrow[dd, "{\cong,k_*}"] 
 & {} \arrow[ldd, "(c)", phantom] 
 & K_0(A) \arrow[rr, "{\cong,k_*}"] \arrow[ll, crossing over, pos=0.3, swap, "{\cong,(i_1)_*}"] 
 & 
 & {K_0(\K(\F_{X_1}))} \arrow[lllldd, "\rho_*"] 
 \\
%
 & 
 & 
 & {} \arrow[ddd, "(e)", phantom] 
 & 
 & {} 
 & 
 & {} \arrow[ddd, "(f)", phantom] 
 & 
 & 
 \\
%
 & 
 & 
 & 
 & 
 & {K_0(\K(\F_{X_1\otimes \TT_{X_2}}))} \arrow[d, "\partial", leftarrow] 
 & 
 & 
 & 
 & 
 \\
%
  & 
  & 
  & 
  & 
  & {K_1(\OO_{X_1\otimes\TT_{X_2}})} \arrow[dd, "{\cong,J_*}"] 
  & 
  & 
  & 
  & 
  \\
%
 & 
 & 
 & {} 
 & 
 & {} 
 & 
 & {} 
 & 
 & 
 \\%
   {K_1(\OO_{X_1})} \arrow[uuuuu, "\partial"] \arrow[rrrrruu, "\theta_*"] \arrow[rrrr, "{\cong,k_*}"'] 
 & 
 & 
 & 
 & {K_1(\K(\F_{X_2\otimes\OO_{X_1}}))} \arrow[r, swap, "\ell_*"] \arrow[ruu, pos=0.45, "(b)", phantom] 
 & {K_1(\TT_{X_2\otimes\OO_{X_1}})} 
 & {} \arrow[luu, "(d)", phantom] 
 & 
 & 
 & {K_1(\OO_{X_1})}, \arrow[uuuuu, "\partial"] \arrow[llll, crossing over, "{\cong,(i_1)_*}"] \arrow[lllluu, "\sigma_*"] 
\end{tikzcd}
};
\end{tikzpicture}
\]
the subdiagrams (a)--(f) commute.
	\end{proposition}
	
	\begin{proof}
		We break this proof down into the six named subdiagrams. We have relegated some of the longer technical  proofs of commutativity to the appendix.

		\emph{Subdiagram (a):} This diagram is induced by the following diagram of $C^*$-algebras which commutes by Lemma~\ref{lem:subgdiagrama} on page~\pageref{lem:subgdiagrama} in Appendix\ref{sec:appendix}:
		\[
		\begin{tikzcd}
			A \arrow[d, "k"] \arrow[r, "k"]         & {\K(\F_{X_2})} \arrow[r, "\ell"] & {\TT_{X_2}} \arrow[d, "k"]            \\
			{\K(\F_{X_1})} \arrow[rr, "\eta"] &                                        & {\K(\F_{X_1\otimes \TT_{X_2}})}.
		\end{tikzcd}\]
		
		 \emph{Subdiagram (b):} This diagram is induced by the following diagram of $C^*$-algebras which commutes by  Lemma~\ref{lem:subgdiagramb} on page~\pageref{lem:subgdiagramb}:

		\[
		\begin{tikzcd}
			{\OO_{X_1}} \arrow[d, "\theta"] \arrow[r, "k"]   & {\K(\F_{X_2\otimes\OO_{X_1}})} \arrow[d, "\ell"] \\
			{\OO_{X_1\otimes\TT_{X_2}}} \arrow[r, "J"] & {\TT_{X_2\otimes\OO_{X_1}}}.
		\end{tikzcd}
		\]

		 \emph{Subdiagram (c):}  This diagram is induced by the following diagram of $C^*$-algebras which commutes by Lemma \ref{lem:subgdiagramc} on page~\pageref{lem:subgdiagramc}:
				\[
		\begin{tikzcd}
			A \arrow[d, "i_1"] \arrow[r, "k"] & {\K(\F_{X_1})} \arrow[d, "\rho"]     \\
			{\TT_{X_2}} \arrow[r, "k"]  & {\K(\F_{X_1\otimes\TT_{X_2}})}.
		\end{tikzcd}
		\]

		 \emph{Subdiagram (d):}  It suffices to show that  the following diagram  of $C^*$-algebras
			\[
		\begin{tikzcd}
			{\OO_{X_1}} \arrow[r, "i_1"] \arrow[rd, "\sigma"] & {\TT_{X_2\otimes\OO_{X_1}}}                \\
			& {\OO_{X_1\otimes\TT_{X_2}}} \arrow[u, "J"]
		\end{tikzcd}
		\]
commutes.
Fix $\xi=\xi'\cdot \langle\xi', \xi'\rangle\in X_1$ and $a\in A$. Since $\sigma$ is induced by $(i_1, \varphi)$, the formulas for $J$ from  Lemma~\ref{lem:fletcherisomorphism} yield
\begin{align*}
J\circ \sigma(j_2(\xi)) &= J(j_2(\varphi(\xi))=J(j_2(\xi'\otimes i_1(\langle\xi', \xi'\rangle)))=i_1(j_2(\xi)) \text{\ and\ }\\
J\circ\sigma(j_1(a)) &= J(j_1(i_1(a)))=i_1(j_1(a)).
\end{align*}
Hence $J\circ \sigma=i_1$.

 \emph{Subdiagram (e):}
Let $(\ell\circ k,\psi)\colon  X_1\to X_1\otimes_A \TT_{X_2}$ be the right-Hilbert $A$-bimodule morphism
from Lemma~\ref{lem:fcbimodulemorphism}.
 Let $\eta\colon \K(\F_{X_1})\to \K(\F_{X_1\otimes\TT_{X_2}})$,  $\theta\colon  \TT_{X_1}\to  \TT_{X_1\otimes\TT_{X_2}} $ and $\Upsilon\colon  \OO_{X_1}\to  \OO_{X_1\otimes\TT_{X_2}}$  be the induced homomorphisms.  Then 		
\[
\begin{tikzcd}
0 \arrow[r] & \K(\F_{X_1}) \arrow[r, "\ell"] \arrow[d,"\eta"] &  \TT_{X_1} \arrow[r, "q"] \arrow[d,"\Upsilon"] &  \OO_{X_1} \arrow[r] \arrow[d, "\theta"] & 0 \\
0 \arrow[r] & \K(\F_{X_1\otimes\TT_{X_2}}) \arrow[r, "\ell"] & \TT_{X_1\otimes\TT_{X_2}} \arrow[r,"q"] & \OO_{X_1\otimes\TT_{X_2}} \arrow[r] & 0
\end{tikzcd}
\]
commutes. It follows from the naturality of the index map as in \cite[Proposition~9.1.5]{RordamLL2000} that
\[
\begin{tikzcd}
K_0(\K(\F_{X_1}))\arrow[r, "\eta"] & K_0( \K(\F_{X_1\otimes\TT_{X_2}})) \\
K_1( \OO_{X_1} ) \arrow[u, "\partial"] \arrow[r, "\theta"] & K_1(\OO_{X_1\otimes\TT_{X_2}}) \arrow[u, "\partial"]
\end{tikzcd}
\]
commutes,  and hence subdiagram (e) commutes as well.

 \emph{Subdiagram (f):}  As for diagram (e), this diagram commutes by naturality
with respect to the right-Hilbert bimodule morphism $\varphi\colon  X_1\to X_1\otimes_A \TT_{X_2} $ from Lemma~\ref{lem:incbimodulemorphism}.
\end{proof}

\begin{lemma}\label{lem:compositionfrom Pimsner}  Let $E$ be a right-Hilbert $A$-bimodule with faithful left action.
 Then \[(\ell\circ k)_*=(i_1)_*\circ (1-[E]).\]
\end{lemma}
\begin{proof}
We show that the desired identity follows from \cite[Lemma~4.7]{Pimsner}.
Most of the work boils down to reconciling Pimsner's notation with ours.

Pimsner's $KK$-class $\alpha \in KK(A, \TT_E)$ is the class of the inclusion $i_1 \colon A \hookrightarrow \TT_E$ (see \cite[Definition~4.1]{Pimsner}), and \cite[Theorem~4.4]{Pimsner} shows that the class $\beta \in KK(\TT_E, A)$ appearing in \cite[Lemma~4.7]{Pimsner} is inverse to $\alpha$. Hence the conclusion of \cite[Lemma~4.7]{Pimsner} can be rearranged as
\begin{equation}\label{eq:rearranged once}
[\mathcal{E}_{+, I}] \hatimes (\iota_I - [E]) \hatimes \alpha = [j].
\end{equation}
Pimsner's $\mathcal{E}_+$ is the Fock space $\F(E)$ (see the first displayed equation on \cite[page~191]{Pimsner}). The ideal $I$ in \cite[Lemma~4.7]{Pimsner} is the preimage of $\K(E)$ under the homomorphism that implements the left action; in our instance this is all of $A$, because we assume a compact left action. Hence $\iota_I$ is the map induced by the identity homomorphism on $A$, which is the identity map $1 = 1_{K_0(A)}$; so $\iota_I - [E]$ is $1 - [E]$. Pimsner's $\mathcal{E}_{+, I}$ is defined on \cite[page~205]{Pimsner} as the set $\{\xi \in \mathcal{E}_+ : \langle \xi, \xi\rangle \in I\}$, so in our instance is equal to $\F(E)$. So the Kasparov class $[\mathcal{E}_{+,I}]$ appearing in \cite[Lemma~4.7]{Pimsner} is the Kasparov class in $KK(\K(\F(E)), A)$ determined by the imprimitivity bimodule $\F(E)$. The conjugate module, which implements the inverse $KK$-class, is the imprimitivity bimodule $\F(E)^*$ determined by the inclusion of $A$ as the full corner of $\K(\F(E))$ corresponding to the summand $A = E^{\otimes 0} \in \F(E)$. This inclusion is the map that we call $k$, and so the Kasparov class $[k] \in KK(A, \F(E))$ is inverse to $[\mathcal{E}_{+, I}]$. So we can rearrange~\eqref{eq:rearranged once} as
\begin{equation}\label{eq:rearranged twice}
(1 - [E]) \hatimes \alpha = [k] \hatimes [j].
\end{equation}
Pimsner writes $j$ for the inclusion of $\K(\F_X)$ in $\TT_X$ that we call $\ell$. By \cite[Proposition~18.7.2(a)]{Blackadar}, we have $[k] \hatimes [\ell] = [\ell \circ k]$.
Since  and $ (i_1)_*(1 - [E])=(1 - [E]) \hatimes [i_1] =(1 - [E]) \hatimes \alpha$, \eqref{eq:rearranged twice} becomes
\begin{equation}\label{eq:rearranged thrice}
(i_1)_*(1 - [E]) = [\ell \circ k],
\end{equation}
as elements of $KK(A, \TT_E)$. In particular, the maps
\[
\cdot \hatimes (i_1)_*(1 - [E]) \colon  KK(\C, A) \to KK(\C, \TT_E)\quad\text{and}\quad
\cdot \hatimes [\ell \circ k] \colon KK(\C, A) \to KK(\C, \TT_E).
\]
coincide. Composing the first of these maps with the isomorphisms $K_*(A) \cong KK^*(\C, A)$ and $K_*(\TT_E) \cong KK^*(\C, \TT_E)$ of \cite[Corollary~18.5.4]{Blackadar} is the definition of the map $(i_1)_*\circ (1 - [E])$ appearing in the statement of the lemma. By \cite[Proposition~18.7.2(a)]{Blackadar} again, the map $\cdot \hatimes [\ell \circ k] \colon  KK(\C, A) \to KK(\C, \TT_E)$ is the map obtained from functoriality of $KK$ described in \cite[Section~17.8]{Blackadar}; there this map is denoted $(\ell \circ k)_*$, but to avoid confusion with the notation in the statement of the lemma, here we will denote it by $(\ell \circ k)\widetilde{\;}$. Given a nondegenerate homomorphism $\phi \colon  B \to C$ of $C^*$-algebras, we have $\ell^2(B) \otimes_B {_\phi C_C} \cong \ell^2(C)_C$ as right Hilbert modules via the isomorphism $(b_i)_{i \in \N} \otimes c = (\phi(b_i)c)_{i \in \N}$, and this isomorphism carries the left action by scalar multiples of a projection $\big(p_{i,j}\big) \in M_n(B)$ to the left action by scalar multiples of $\big(\phi(p_{i,j})\big)$. Hence the description in  \cite[Section~17.8]{Blackadar} shows that the isomorphisms $K_*(A) \cong KK^*(\C, A)$ and $K_*(\TT_E) \cong KK^*(\C, \TT_E)$ of \cite[Corollary~18.5.4]{Blackadar} carry $(\ell \circ k)\widetilde{\;}$ to the homomorphism in $K$-theory induced by $\ell \circ k$, namely $(\ell \circ k)_*$. So, after identifying $K$-theory with $KK(\C, \cdot)$ as in  \cite[Corollary~18.5.4]{Blackadar}, Equation~\ref{eq:rearranged thrice} becomes the desired equation $(\ell\circ k)_*=(i_1)_* \circ (1-[E])$.
\end{proof}
		
	\begin{proof}[Proof of Theorem \ref{thm:diagram1}]
	 The KK-diagram behind diagram~\eqref{eq:monster diagram}  is the following, where all left actions are by compacts:
		\begin{equation}\label{eq:KKdiagram}
		\begin{tikzpicture}[yscale=0.85,xscale=0.8,baseline=(a), every node/.style={scale=0.9}]
			\node (a) at (0,0){
				\begin{tikzcd}[row sep=small, column sep = small]
					&     &  {KK_1(\C,\OO_{Y_1})} \arrow[ddddll, bend right=25, "\partial", out=315]  \arrow[ddddrrrrrr, bend left=20, leftarrow, pos=0.3, swap, "\cdot \hatimes \OO_{Y_1}"]   &  &                            &  & {KK_0(\C,A)} \arrow[rr, "1- \cdot\hatimes X_2"]                                        &     & {KK_0(\C,A)}                                        \\
					&     &                                                        &  &                                             &  &                                                                                            & (4) &                                                     \\
					0                                                                                             &     & 0                                                      &  &                                             &  & {KK_0(\C,A)} \arrow[rr, crossing over, "1- \cdot\hatimes X_2"] \arrow[uu, crossing over, pos=0.4, "1- \cdot\hatimes X_1"] &     & {KK_0(\C,A)} \arrow[uu, pos=0.4, "1- \cdot\hatimes X_1"] \\
					&     &                                                        &  &                                             &  &                                                                                            & (3) &                                                     \\
					{KK_0(\C,\OO_{X_1})} \arrow[uu, shift left] \arrow[rr, "1-\cdot\hatimes Y_1"]                  &     & {KK_0(\C,\OO_{X_1})} \arrow[uu] \arrow[rr, "\cdot \hatimes \OO_{Y_1}"]  &  & {KK_0(\C,\OO_{Y_1})} \arrow[rr, "\partial"] &  & {KK_1(\C,\OO_{X_1})} \arrow[uu, "\partial"] \arrow[rr, "1-\cdot\hatimes Y_1"]               &     & {KK_1(\C,\OO_{X_1})}. \arrow[uu, "\partial"]        \\
					& (2) &                                                        &  &                                             &  &                                                                                            &     &                                                     \\
					{KK_0(\C,A)} \arrow[uu, "\cdot\hatimes \OO_{X_1}"] \arrow[rr, "1- \cdot\hatimes X_2"] &     & {KK_0(\C,A)} \arrow[uu, "\cdot\hatimes \OO_{X_1}"] &  &                                             &  & 0 \arrow[uu]                                                                               &     & 0 \arrow[uu]                                        \\
					& (1) &                                                        &  &                                             &  &                                                                                            &     &                                                     \\
					{KK_0(\C,A)} \arrow[rr, "1- \cdot\hatimes X_2"] \arrow[uu, "1- \cdot\hatimes X_1"]    &     & {KK_0(\C,A)} \arrow[uu, "1- \cdot\hatimes X_1"]    &  &                                             &  &                                                                                            &     &
				\end{tikzcd}
			};
		\end{tikzpicture}
		\end{equation}
We will show  diagram~\eqref{eq:KKdiagram} commutes, which implies that~\eqref{eq:monster diagram} commutes. 	
The middle row, including the curved arrow from extreme right to extreme left, is the six-term Pimsner sequence. 
Since $K_1(A)=0$, the  Pimsner sequence associated to the right-Hilbert bimodule $X_1$ over $A$ is
\[
\begin{tikzcd}
			0 \arrow[r] & {KK_1(\C,\OO_{X_1})} \arrow[r, "\partial"] & {KK_0(\C,A)} \arrow[r, "{1-\cdot\hatimes X_1}"] & [2em] {KK_0(\C,A)} \arrow[r, "\cdot\hatimes \OO_{Y_1}"] & {KK_0(\C,\OO_{X_1})} \arrow[r] & 0.
		\end{tikzcd}
		\]
The left two columns in diagram~\eqref{eq:KKdiagram} consist of the last three terms and the right two columns consist of the first three terms of this sequence. It remains to prove that squares (1)--(4) commute. 
		
		 \emph{Squares (1)~and~(4):} These commute because $X_1\otimes_A X_2 \cong X_2 \otimes_A X_1$.
		
		 \emph{Square~(2):} The identity maps, denoted $1$, clearly commute with ${\cdot \hatimes \OO_{X_1}}$. Tracing Square~(2) along the bottom and right gives
		 \[\cdot  \hatimes  (X_2 \hatimes \OO_{X_1}) = \cdot \hatimes  (X_2\otimes_A \OO_{X_1}),\] whereas tracing along the left and top gives \[ \cdot \hatimes (\OO_{X_1}\hatimes  Y_1) = \cdot \hatimes (\OO_{X_1}\hatimes  (X_2\otimes_A \OO_{X_1})) = \cdot\hatimes  (\OO_{X_1}\otimes_{\OO_{X_1}}(X_2\otimes_A \OO_{X_1})).\]
		Because all left actions involved are compact, the associated Fredholm operator can be chosen to be $0$, and it suffices to show that
		\begin{equation}\label{square2}
			X_2 \otimes_A \OO_{X_1} \cong \OO_{X_1}\otimes_{\OO_{X_1}}(X_2\otimes_A \OO_{X_1})\end{equation}
		as right-Hilbert $A$--$\OO_{X_1}$-bimodules. Because $Y_1$ is a right-Hilbert bimodule over $\OO_{X_1}$, $\OO_{X_1}\otimes_{\OO_{X_1}} Y_1 \cong Y_1$ and thus \eqref{square2} holds as an isomorphism of right-Hilbert $\OO_{X_1}$-modules. But because the left action on $X_2\otimes_A \OO_{X_1}$ restricts to multiplication by $j_1(A)$, \eqref{square2} is an isomorphism of right-Hilbert $A$--$\OO_{X_1}$-bimodules as well. 
		
		 \emph{Square (3):}
Since squares (a)--(f) in the diagram of Proposition~\ref{prop:square3} commute, the whole digram commutes.  By two applications of Lemma~\ref{lem:compositionfrom Pimsner},  the compositions $(i_1)^{-1}_*\circ \ell*\circ k_*$ along the top and bottom of the rectangle are $1-[X_2]$  and $1-[Y_1]$, respectively.  After identifying $\delta$ with $k\circ\delta$, it follows that  the square in \eqref{eq:monster diagram} corresponding to Square~(3) commutes (and hence (3) commutes). It follows that diagram~\eqref{eq:KKdiagram} commutes, and hence  diagram~\eqref{eq:monster diagram} commutes as well. 	

It remains to show that diagram~\eqref{eq:monster diagram} is natural with respect to the right-Hilbert bimodule morphisms of Lemma~\ref{lem:fletchermorphisms}; again it suffices to show that diagram~\eqref{eq:KKdiagram} is natural.
Pimsner's six-term sequence in $KK$-theory is natural for covariant right-Hilbert bimodule morphisms by \cite[Propositions 6.4.1~and~7.0.4]{Hume}.

Since the central six-term sequence is the Pimsner sequence associated to the right-Hilbert bimodule $Y_1$, it is natural  with respect to the second pair of homomorphisms in Lemma~\ref{lem:fletchermorphisms}; thus in particular the central horizontal five-term  sequence is natural. Similarly, the four vertical sequences are natural with respect to the first pair in Lemma \ref{lem:fletchermorphisms} because they come from independent six-term Pimsner sequences.
We show that the square
		\begin{equation}\label{eq:square5}
		\begin{tikzcd}[column sep = large]
			{KK_0(\C,A)} \arrow[r, "1-\cdot \hatimes X_2"]                                                                      & {KK_0(\C,A)}                                   \\
			{KK_0(\C,I)} \arrow[r, "1-\cdot \hatimes (X_2 I)"] \arrow[u, "\cdot\hatimes A"] \arrow[ru, "", phantom] & {KK_0(\C,I)} \arrow[u, "\cdot\hatimes A"]
		\end{tikzcd}
		\end{equation}		
		commutes. In order to show that \eqref{eq:square5} commutes, it suffices to show that the essential submodules $I\cdot {}_{I}(A\otimes_A X_2)_A$  and $I\cdot {}_{I}(X_2 I\otimes_I A)_A$ of ${}_{I}(A\otimes_A X_2)_A$ and ${}_{I}(X_2 I\otimes_I A)_A$ are isomorphic. Now ${}_A(A\otimes_A X_2)_A \cong {}_A (X_2)_{A}$ as right-Hilbert bimodules, but because $I$ is invariant, $I X_2 \subset X_2 I$ and hence the essential submodule of ${}_I(X_2)_A$ is equal to ${}_I (X_2)_{I}$. Finally \[{}_{I}(X_2 I\otimes_I A)_A \cong {}_I(X_2 \otimes_I I)_A \cong {}_I(X_2 \otimes_I I)_I \cong {}_I(X_2)_I\] and we have shown that~\eqref{eq:square5} commutes.  The square
		\[
		\begin{tikzcd}
			{KK_0(\C,\OO_{X_1})} \arrow[r, "1-\cdot\hatimes Y_1"]                                                                              & {KK_0(\C,\OO_{X_1})}                                                  \\
			{KK_0(\C,\OO_{X_1 I})} \arrow[u, "\cdot\hatimes \OO_{X_1}"] \arrow[r, "1-\cdot\hatimes Y_1^I"] \arrow[ru, "", phantom] & {KK_0(\C,\OO_{X_1I}),} \arrow[u, "\cdot\hatimes \OO_{X_1}"]
		\end{tikzcd}\] commutes by naturality of the central horizontal five-term sequence. Thus the entire diagram is natural with respect to the morphisms from Lemma \ref{lem:fletchermorphisms}.
	\end{proof}

	
	Let $A$ be  a $C^*$-algebra, and let $X_1$ and $X_2$ be Hilbert $A$-bimodules  with faithful compact left actions such that $X_1\otimes_A X_2 \cong X_2\otimes_A X_1$.
	Consider the homomorphism $1- [X_2] \colon  K_0(A) \to K_0(A)$.   Since $[X_1] = \cdot \hatimes X_1$ and $[X_2]= \cdot \hatimes X_2$ commute, $1- [X_2]$ descends to a homomorphism
	\[1-\widetilde{[X_2]} \colon \coker(1-[X_1])\to\coker(1-[X_1]),\]
    and restricts to a homomorphism
    \[1-[X_2]|_{\ker(1 - [X_1])}\colon \ker(1-[X_1])\to\ker(1-[X_1]).\]
    We show that we can, as in \cite{aHNS},
	rewrite the cokernel and kernel of these two maps respectively in terms of the following maps:
	\[(1-[X_1],1-[X_2]) \colon  K_0(A)\oplus K_0(A) \to K_0(A), \quad (p,q) \mapsto (1-[X_1])p + (1-[X_2])q,\] and
	\[\begin{pmatrix}
		1-[X_1] \\ 1-[X_2]
	\end{pmatrix} \colon  K_0(A) \to K_0(A)\oplus K_0(A), \quad p \mapsto \begin{pmatrix} p-[X_1]p \\ p - [X_2]p\end{pmatrix}.\]

	 \begin{lemma}\label{lem:blockmatrixmaps}
	 	The map $p + \image(1-[X_1])  \mapsto p + \image(1-[X_1],1-[X_2])$ gives an isomorphism
	 	\[\coker\bigg(1-\widetilde{[X_2]} \colon \coker(1-[X_1]) \to \coker(1-[X_1])\bigg)  \cong \coker(1-[X_1],1-[X_2]),
	 	\]
	 	and
	 	\[\ker\bigg(1-[X_2]|_{\ker(1 - [X_1])} \colon  \ker(1-[X_1])\to \ker(1-[X_1])\bigg) = \ker\begin{pmatrix}
	 		1-[X_1] \\ 1- [X_2]
	 	\end{pmatrix}.\]
	 \end{lemma}
	
	 \begin{proof}
	 	For $p\in K_0(A)$,
	 	\begin{align*}
	 		p+ \image(1-[X_1]) & \in \image(1-\widetilde{[X_2]}) \\
	 		& \iff p - (1-[X_1])q \in \image(1-[X_2]) \text{ for some } q\in K_0(A) \\
	 		& \iff p = (1-[X_1])q + (1-[X_2])t \text{ for some } p,q,t \in K_0(A) \\
	 		& \iff p \in \image(1-[X_1],1-[X_2]).
	 	\end{align*}
	 	Thus, $p + \image(1-[X_1])) \mapsto p + \image(1-[X_1],1-[X_2])$, which is clearly surjective from $\coker(1-\widetilde{[X_2]})$ to $\coker(1-[X_1],1-[X_2])$, is also injective, and hence is an isomorphism. Since,
	 	\[\ker(1-[X_2]|_{\ker(1 - [X_1])}) = \ker(1-[X_1])\cap\ker(1-[X_2]) = \ker\begin{pmatrix}
	 		1-[X_1] \\ 1- [X_2]
	 	\end{pmatrix},\] we are done.
	 \end{proof}
	
	We now have the tools to prove our main $K$-theoretic result which is a generalisation of \cite[Theorem 4.20]{aHNS} for general commuting right-Hilbert bimodules with compact left actions and coefficient algebra with trivial $K_1$-group.  Equation~\ref{eq:thmK_0} generalises an anologous exact sequence proved in \cite[Theorem~6.10]{FarsiKumjianPaskSims} for rank-$2$ Deaconu--Renault groupoids, but~\eqref{eq:thm_idealinclusion} is new even for  rank-$2$ Deaconu--Renault groupoids---see Theorem~\ref{thm:DRkgroups2} below.

	\begin{theorem}\label{thm:naturalses} Let $A$ be a separable $C^*$-algebra with $K_1(A)=0$.
		Let $X_1$ and $X_2$ be countably generated right-Hilbert $A$-bimodules with faithful compact left actions of $A$ such that $X_1\otimes_A X_2 \cong X_2\otimes_A X_1$. Let $Y_1$ denote the right-Hilbert bimodule ${_{\OO_{X_1}}\!(X_2 \otimes_A \OO_{X_1})_{\OO_{X_1}}}$ from the Deaconu--Fletcher construction of \S\ref{sec:fletcher}, and
write $j_1^{(1)} \colon  A \to \OO_{X_1}$ and $j_1^{(2)} \colon  \OO_{X_1}\to \OO_{Y_1}$ for the canonical inclusion maps, and define $j\colon  \coker(1-[X_1],1-[X_2]) \to K_0(\OO_{Y_1})$ by
		\[
		j(p+ \image(1-[X_1],1-[X_2])) = [(j_1^{(2)}\circ j_1^{(1)})_*(p)],\quad\text{for all $p \in K_0(A)$.}
		\]
		Then the composition $\tau \coloneqq \partial \circ \partial \colon  K_0(\OO_{Y_1}) \to K_0(A)$ in the diagram~\eqref{eq:monster diagram} makes the sequence%
		\begin{equation}\label{eq:thmK_0}
		\begin{tikzcd}
		0 \arrow[r] & {\coker(1-[X_1],1-[X_2])} \arrow[r, "j"] & K_0(\OO_{Y_1}) \arrow[r, "\tau"] & {\ker{\begin{pmatrix}1-[X_1]\\ 1-[X_2]\end{pmatrix}}} \arrow[r] & 0
		\end{tikzcd}
		\end{equation}
		exact.
		Suppose that $I$ is an $X_1$- and $X_2$-invariant ideal and let $Y_1^I$ denote the right-Hilbert bimodule ${}_{\OO_{X_1 I}}(X_2 I\otimes_I \OO_{X_1 I})_{\OO_{X_1 I}}$. Let $\imath \colon  I \hookrightarrow A$ be the inclusion. Then there is a map
		\[\tilde{\imath}\colon  \coker(1-[X_1 I],1-[X_2 I]) \to \coker(1-[X_1],1-[X_2]),\]
		such that
		\[\tilde{\imath}(p + \image(1-[X_1 I],1-[X_2 I])) = \imath_*(p) + \image(1-[X_1],1-[X_2]).\]
		Let $\phi \colon  \OO_{Y_1^I}\to \OO_{Y_1}$ be the homomorphism of Lemma~\ref{lem:fletchermorphisms} and $j_I$ be the map analogous to $j$ from $\coker(1-[X_1I], 1-[X_2I])\to K_0(\OO_{Y_1^I})$. Then
		the diagram
		\begin{equation}\label{eq:thm_idealinclusion}
		\begin{tikzcd}
		0 \arrow[r] & {\coker(1-[X_1 I],1-[X_2 I])} \arrow[r, "j_I"] \arrow[d, "\tilde{\imath}"] & K_0(\OO_{Y_1^I}) \arrow[r, "\tau_I"] \arrow[d, "\phi_*"] & {\ker{\begin{pmatrix}1-[X_1 I]\\ 1-[X_2 I]\end{pmatrix}}} \arrow[r] \arrow[d, "\imath_*"] & 0 \\
		0 \arrow[r] & {\coker(1-[X_1],1-[X_2])} \arrow[r, "j"]                                             & K_0(\OO_{Y_1}) \arrow[r, "\tau"]                       & {\ker{\begin{pmatrix}1-[X_1]\\ 1-[X_2]\end{pmatrix}}} \arrow[r]                                   & 0
		\end{tikzcd}
		\end{equation}
		commutes and has exact rows.
	\end{theorem}
	\begin{proof}
		Consider the left part of the commuting diagram in Theorem \ref{thm:diagram1}:
		\[
		\begin{tikzcd}
			K_0(A) \arrow[d, "{1-[X_2]}"] \arrow[r, "{1-[X_1]}"] & K_0(A) \arrow[d, "{1-[X_2]}"] \arrow[r, "(j_1^{(1)})_*"] & {K_0(\OO_{X_1})} \arrow[d, "{1-[Y_1]}"] \arrow[r] & 0 \\
			K_0(A) \arrow[r, "{1-[X_1]}"]                              & K_0(A) \arrow[r, "(j_1^{(1)})_*"]                              & {K_0(\OO_{X_1})} \arrow[r]                        & 0.
		\end{tikzcd}\]
    Let $(j_1^{(1)})^{\!\thicksim}_* \colon  \coker(1-[X_1]) \to K_0(\OO_{X_1})$ be the isomorphism induced by $(j_1^{(1)})_*$. By \cite[Remark~2.2]{aHNS} we obtain the commuting square
		\[
		\begin{tikzcd}
			{\coker(1-[X_1])} \arrow[d, "{1-\widetilde{[X_2]}}"] \arrow[r, "(j_1^{(1)})^{\!\thicksim}_*{,}\cong"] & [2em] {K_0(\OO_{X_1})} \arrow[d, "{1-[Y_1]}"] \\
			{\coker(1-[X_1])} \arrow[r, "(j_1^{(1)})^{\!\thicksim}_*{,}\cong"]                                & {K_0(\OO_{X_1}).}
		\end{tikzcd}\]
		Likewise, let $\widetilde{\partial} : K_1(\OO_{X_1}) \to \ker(1-[X_1])$ be the isomorphism induced by $\partial$. Then the right part
		\[
		\begin{tikzcd}
			0 \arrow[r] & {K_1(\OO_{X_1})} \arrow[r, "\partial"] \arrow[d, "{1-[Y_1]}"] & K_0(A) \arrow[r, "{1-[X_1]}"] \arrow[d, "{1-[X_2]}"] & K_0(A) \arrow[d, "{1-[X_2]}"] \\
			0 \arrow[r] & {K_1(\OO_{X_1})} \arrow[r, "\partial"]                        & K_0(A) \arrow[r, "{1-[X_1]}"]                              & K_0(A)
		\end{tikzcd}
		\]
		 of the commuting diagram in Theorem \ref{thm:diagram1} shows that the diagram
		\[
		\begin{tikzcd}
			{K_1(\OO_{X_1})} \arrow[d, "{1-[Y_1]}"] \arrow[r, "\cong"] & {\ker(1-[X_1])} \arrow[d, "{1-[X_2]|_{\ker(1-[X_1])}}"] \\
			{K_1(\OO_{X_1})} \arrow[r, "\cong"]                        & {\ker(1-[X_1])}
		\end{tikzcd}
		\]
		commutes. Both of these new squares are natural with respect to the desired right-Hilbert bimodule morphisms.
		Hence we obtain the commuting diagram
		\[
		\begin{tikzpicture}[scale=0.95,baseline=(a), every node/.style={scale=0.95}]
			\node (a) at (0,0){
				\begin{tikzcd}
					&                                             &                                                                    & {\ker(1-[X_1])} \arrow[r, "{1-[X_2]|_{\ker(1-[X_1])}}"]              &[+3.2em] {\ker(1-[X_1])}                     \\
					{K_0(\OO_{X_1})} \arrow[r, "{1-[Y_1]}"]                                       & {K_0(\OO_{X_1})} \arrow[r, "(j_1^{(2)})_*"] & {K_0(\OO_{Y_1})} \arrow[r, "\partial"] & {K_1(\OO_{X_1})} \arrow[u, "\cong"] \arrow[r, "{1-[Y_1]}"] & {K_1(\OO_{X_1})}. \arrow[u, "\cong"] \\
					{\coker(1-[X_1])} \arrow[u, "(j_1^{(1)})_*{,}\cong"] \arrow[r, "{1-\widetilde{[X_2]}}"] & {\coker(1-[X_1]).} \arrow[u, "(j_1^{(1)})_*{,}\cong"]  &                                                                    &                                                                  &
				\end{tikzcd}
			};
		\end{tikzpicture}
		\]
		Let $\jmath = (j_1^{(2)})_*\circ (j_1^{(1)})^{\!\thicksim}_*$. Then we get the exact sequence
		\[
		\begin{tikzpicture}[scale=0.95,baseline=(a), every node/.style={scale=0.95}]
			\node (a) at (0,0){
				\begin{tikzcd}
					{\coker(1-[X_1])} \arrow[r, "{1-\widetilde{[X_2]}}"] & [1em] {\coker(1-[X_1])} \arrow[r,"\jmath"] &[-1em] {K_0(\OO_{Y_1})} \arrow[r] &[-1em] {\ker(1-[X_1])} \arrow[r, "{1-[X_2]|_{\ker(1-[X_1])}}"] &[3.2em] {\ker(1-[X_1])}.
				\end{tikzcd}
			};
		\end{tikzpicture}
		\]
		Let $j : \coker(1-\widetilde{[X_2]}) \to K_0(\OO_{Y_1})$ be the homomorphism induced by $\jmath$. Applying the first isomorphism theorem to both sides of the above sequence gives the short exact sequence
				\[
				\begin{tikzcd}
						0 \arrow[r] & {\coker(1-\widetilde{[X_2]})} \arrow[r, "j"] & K_0(\OO_{Y_1}) \arrow[r] & {\ker(1-[X_2]|_{_{\ker(1-[X_1])}})} \arrow[r] & 0.
					\end{tikzcd}
				\]
		By Lemma \ref{lem:blockmatrixmaps}, composing by the appropriate isomorphisms and relabelling $j$ gives the desired short exact sequence and the description of $j$. Since Square~(5) in the proof of Theorem~\ref{thm:diagram1} commutes, $\tilde{\imath}\colon \coker(1-[X_1I],1-[X_2 I])\to \coker(1-[X_1],1-[X_2])$ (see Theorem~\ref{thm:naturalses}) is well defined and $\imath_*\colon K_0(I)\to K_0(A)$ maps
        $\ker\big(\begin{smallmatrix} 1-[X_1 I]\\1-[X_2 I]\end{smallmatrix}\big)$ into $\big(\ker\begin{smallmatrix}	1-[X_1]\\1-[X_2] \end{smallmatrix}\big)$. The naturality of the big diagram in Theorem~\ref{thm:diagram1} with respect to the right-Hilbert bimodule morphisms of Lemma \ref{lem:fletchermorphisms} completes the proof.
	\end{proof}
	
\section{\texorpdfstring{$K$}{K}-theory for continuous functions on totally disconnected spaces}\label{sec:KK}

The purpose of this self-contained section is to prove Theorem~\ref{thm:tdlcKtheory}, which we need to pull over the abstract $K$-theoretic results from \S\ref{sec:ktheoryfletcher} to the Hilbert modules associated to a rank-2 Deaconu--Renault groupoid.

Throughout this section, let $X$ and $Y$ be  second-countable totally disconnected locally compact Hausdorff spaces, and let $T\colon X\to Y$ be a surjective local homeomorphism.
We now define the maps appearing in the commuting diagram of Theorem~\ref{thm:tdlcKtheory}, starting with the right-Hilbert $C_0(Y)$-module $E_T$.
For $ \xi, \eta \in C_c(X)$ and $b\in C_0(Y)$ define
\[
\langle \xi, \eta  \rangle_{C_0(Y)}(y) = \sum_{Tx = y} \overline{\xi(x)}\eta(x)\text{\ and \ } (\xi\cdot b)(x) = \xi(x)b(Tx).
\]
Let $E_T$ be the  completion of $C_c(X)$ under the norm induced by the inner product $\langle \cdot, \cdot  \rangle_{C_0(Y)}$.  Then $E_T$ is a full right-Hilbert  $C_0(Y)$-module.  In addition, for $ \xi  \in C_c(X)$ and $a\in C_0(X)$ define   $a\cdot\xi(x) = a(x)\xi(x)$; this extends to an action of $C_0(X)$ on $E_T$ by compact operators.   To see the last statement, let $U$ be a  compact and open subset of $X$ such that $T|_U$ is injective.  Then
\begin{align}\label{action by compacts}
\big(\Theta_{1_U, 1_U}(\xi)\big)(x)
    &= \big(1_U\cdot\langle 1_U, \xi\rangle_{C_0(Y)}\big)(x)
     = 1_U(x)\sum_{T(z)=T(x)} \overline{1_U(z)}\xi(z)\\
    &= 1_U(x)\xi(x)=(1_U\cdot\xi)(x)\notag.
\end{align}

By \cite[Exercise~3.4]{RordamLL2000}, for each compact open subset $V
\subset X$, there is an isomorphism from $C(V, \Z)$ to $K_0(C(V))$ that
carries the indicator function $1_U$ of a compact open set $U \subset V$ to
the $K_0$ class $[1_U]_0$. Continuity of $K$-theory implies that there is an isomorphism $\sigma^X \colon  C_c(X, \Z) \to K_0(C_0(X))$
such that
\begin{equation}\label{eq:sigma def}
	\sigma^X(1_U) = [1_U]_0\text{ for each compact open } U \subset X.
\end{equation}
By Theorem~3.8 of \cite{BKR} there is an isomorphism
\[
\Theta^X\colon  K_0(C_0(X)) \to KK_0(\C, C_0(X))
\]
such that, for   each
compact open set $U \subset \Omega$,
\begin{equation*}
\Theta^X([1_U]_0)= \big[z \mapsto z1_{C(U)}, C(U) \oplus 0, 0,
\id\oplus(-\id)\big];
\end{equation*}
the latter is the class of the Kasparov module obtained by regarding $C(U)$ as a trivially graded right-Hilbert module over itself with left action of $z\in\C$ as multiplication by $z1_U$.  (We refer the reader to \cite{Blackadar} for details on Kasparov modules and their product.)
Let $T_*\colon  C_c(X,\Z) \to C_c(Y,\Z)$ be the map
\begin{equation}\label{eq:T* def}
    (T_* h)(y) = \sum_{Tx = y} h(x);
\end{equation}
if $U$ is a compact open set for which $T|_U  \colon  U \to T(U)$ is a homeomorphism, then $T_*(1_U) = 1_{T(U)}$.

%

\begin{theorem}\label{thm:tdlcKtheory}
Let $X,Y$ be second-countable totally disconnected locally compact Hausdorff spaces, and let $T\colon X\to Y$ be a surjective local homeomorphism. Let $\sigma^X$, $\Theta^X$, $T_*$ and $\cdot\hatimes E_T$ be as described above, and let $[E_T] : K_0(C_0(X)) \to K_0(C_0(Y))$ be the map of~\eqref{eq:KKnotation}. Then
the diagram
\[
 \begin{tikzcd}
	{C_c(X,\Z)} \arrow[rrr, "T_*"] \arrow[d, "\sigma^X"'] &
	& &{C_c(Y,\Z)} \arrow[d, "\sigma^Y"]
	\\
	K_0(C_0(X)) \arrow[r, "{\Theta^X}"'] \arrow[rrr,  "{[E_T]}", bend right]
        & {KK_0(\C,C_0(X))} \arrow[r, "\cdot\hatimes E_T"] & {KK_0(\C,C_0(Y))} & K_0(C_0(Y)) \arrow[l, "{\Theta^Y}"]
\end{tikzcd}
\]
commutes, and $[E_T]$ is the unique homomorphism such that $[E_T]([1_U]_0) = [1_{T(U)}]_0$ for every compact open $U$ such that $T|_U$ is a homeomorphism.
\end{theorem}
\begin{proof}
%
The formulas for $\sigma^Y$, $\sigma^X$ and $T_*$ imply that $\sigma^Y \circ T_* \circ (\sigma^X)^{-1}([1_U]_0) = [1_{T(U)}]_0$ whenever $U \subset X$ is compact open and $T|_U$ is a homeomorphism. So it suffices to show that the diagram commutes. That the bottom part commutes is just the definition of $[E_T]$.

Since $E_T$ is a full right-Hilbert $C_0(Y)$-module, we can view it as a  $\K(E_T)$--$C_0(Y)$-imprimitivity bimodule (denoted ${_\mathcal{K}E_T}$) by \cite[Proposition~3.8]{RaeburnWilliams}.  The algebras $\K(E_T)$ and $C_0(Y)$ are complementary full  corners in the linking algebra $L(E_T)$ of $E_T$; we write $P_{\mathcal{K}}, P_Y$ for the multiplier projections onto these corners. By \cite[Proposition~1.2]{Paschke}, the corner inclusions induce isomorphisms
\begin{equation}\label{eq:map names}
    i_{\K(E_T)}\colon K_0(\K(E_T))\to K_0(L(E_T)) \qquad\text{and}\qquad i_{C_0(Y)}\colon K_0(C_0(Y))\to K_0(L(E_T)).
\end{equation}

Fix a  compact open $U\subset X$ such that $T|_U$ is a homeomorphism.  By \eqref{action by compacts}, the homomorphism $\varphi\colon C_0(X) \to \mathcal{L}(E_T)$ implementing the left action satisfies $\varphi(1_U)=\Theta_{1_U, 1_U}$. Let $\varphi_*\colon K_0(C_0(X))\to K_0(\K(E_T))$ be the induced homomorphism. The inclusions of  $\Theta_{1_U, 1_U}={}_ {\K(E_T)}\langle 1_U, 1_U\rangle$ and $1_{T(U)}=\langle 1_U, 1_U\rangle_{C_0(Y)}$ in $L(E_T)$ are Murray--von Neumann equivalent via the partial isometry $\big(\begin{smallmatrix}0&1_U\\0&0\end{smallmatrix}\big)$.
Thus
\[
    i_{C_0(Y)}^{-1}\circ i_{\K(E_T)}\circ\varphi_*([1_U]_0) = i_{C_0(Y)}^{-1}\circ i_{\K(E_T)}([\Theta_{1_U, 1_U}]_0)=[1_{T(U)}]_0.
\]
So it suffices to show that
\begin{equation}\label{eq:required K-th equality}
\Theta_Y^{-1} \circ (\cdot \hatimes E_T) \circ \Theta_X = i_{C_0(Y)}^{-1} \circ i_{\mathcal{K}(E_T)} \circ \varphi_*;
\end{equation}
the final statement then follows by definition of $[E_T]$---see~\eqref{eq:KKnotation}.

We write ${_\varphi \mathcal{K}(E_T)}$ for $\mathcal{K}(E_T)$ regarded as a right-Hilbert $\mathcal{K}(E_T)$-module with the canonical right action and inner-product, and compact injective left action implemented by $\varphi$. Note that $P_{\mathcal{K}} L(E_T)$ and $L(E_T)P_Y$ are Morita equivalences from $\mathcal{K}(E_T)$ to $L(E_T)$ and from $L(E_T)$ to $C_0(Y)$ respectively. We have $P_{\mathcal{K}} L(E_T) = {_{i_{\mathcal{K}(E_T)}} L(E_T)}$ and $L(E_T) P_Y = E_T \oplus C_0(Y)$, which is the conjugate of the module ${_{i_{C_0(Y)}} L(E_T)}$. Since
\begin{align*}
{_\varphi \mathcal{K}(E_T)} \otimes_{\mathcal{K}(E_T)} P_{\mathcal{K}} L(E_T) \otimes_{L(E_T)} L(E_T)P_Y
    &\cong {_\varphi \mathcal{K}(E_T)} \otimes_{\mathcal{K}(E_T)} P_{\mathcal{K}} L(E_T) P_Y\\
    &\cong {_\varphi \mathcal{K}(E_T)} \otimes_{\mathcal{K}(E_T)} {_{\mathcal{K}}E_T}
     \cong E_T,
\end{align*}
the map $\cdot \hatimes E_T$ coincides with $[i_{C_0(Y)}]^{-1} \hatimes [i_{\mathcal{K}(E_T)}] \hatimes \varphi_*$. Conjugating this with the isomorphisms $\Theta^X$ and $\Theta^Y$ of \cite[Theorem~3.8]{BKR} described above gives~\eqref{eq:required K-th equality}.
\end{proof}

\section{\texorpdfstring{Example: $K$}{K}-theory for rank-two Deaconu--Renault groupoids}\label{sec:ktheoryDR}

We will say that a groupoid $\mathcal{G}$ is \emph{ample} if it is \'etale and totally disconnected. In this section we apply the $K$-theoretic results of \S\ref{sec:ktheoryfletcher} and \S\ref{sec:KK} to the $C^*$-algebra of an ample rank-$2$ Deaconu--Renault groupoid and an ideal arising from an invariant open subset of the unit space of the groupoid.  The main result of the section is Theorem~\ref{thm:DRkgroups2}.

For the genesis of groupoid $C^*$-algebras, we refer the reader to \cite{Renault80}. The first examples of what we now call Deaconu--Renault groupoids appeared there \cite[Definition~III.2.1]{Renault80} as models for the Cuntz algebras, and were generalised to graph groupoids in \cite{KPRR}. The general construction, for a single local homeomorphism of a locally compact Hausdorff space, was introduced in \cite{Deaconu-groupoid}. Here we follow the conventions and notation of \cite{SSW} (for Deaconu--Renault groupoids specifically, see \cite[Examples 8.1.16~and~8.3.7]{SSW}).

Throughout this section, let $\Omega$ be a second-countable totally disconnected locally compact Hausdorff space, and let $T_1$ and $T_2$  be commuting surjective local homeomorphisms  on $\Omega$.
Equivalently, $T_1$ and $T_2$ give an action of $\N^{2}$ on $\Omega$ by local homemomorphisms.

\subsection*{Background on Deaconu--Renault  groupoids} For $m=(m_1, m_2)\in\N^2$ we write $T^m$ for $(T_1)^{m_1}(T_2)^{m_2}$. The corresponding \emph{Deaconu--Renault  groupoid} from \cite{Deaconu-groupoid} is the set
  \begin{equation}\label{eq-DRgroupoid}
G_{T}\coloneqq \bigcup_{m,n\in\N^{2}}\{(x,m-n,y) \in \Omega \times \Z^2 \times \Omega : T^{m}x = T^{n}y
\}
  \end{equation}
with unit space $G_{T}^{(0)} = \{(x,0,x) : x \in \Omega\}$ identified with
  $\Omega$, range and source maps $r(x,n,y) = x$ and $s(x,n,y) = y$, and
  operations $(x,n,y)(y,m,z)=(x,n+m,z)$ and $(x,n,y)^{-1}= (y,-n,x)$.
 For compact open sets $U,V \subset \Omega$ and for $m,n \in \N^2$,  set
  \begin{equation}\label{eq:basicset}
    Z(U,m,n,V)\coloneqq
    \{(x,m-n,y):\text{$x\in U$, $y\in V$ and
        $T^{m}x=T^{n}y$}\}.
  \end{equation}
The sets at \eqref{eq:basicset} form a  basis of compact open sets for a locally compact Hausdorff topology on $G_T$. Furthermore, the sets  $Z(U,m,n,V)$ such that in addition $T^m|_U$ and $T^n|_V$ are homemomorphisms and $T^m(U)=T^n(V)$ are a basis for the same topology. With respect to this basis,  $G_T$ is an ample second-countable locally compact Hausdorff groupoid---see, for example, \cite[Lemma~3.1]{SimsWilliams2016} (there $\Omega$ may not be totally disconnected, and so their $G_T$  may only  be  \'etale).

\subsection*{The product system  \texorpdfstring{from \cite[\S5]{CLSV}}{of Carlsen-Larsen-Sims-Vittadello}}	 A topological $2$-graph $\Lambda_T$  can be thought of as either a higher-rank generalisation of a topological graph, or as a topological generalisation of a discrete higher-rank graph; see \cite{Yeendtopk} for details. We define a topological $2$-graph, whose infinite-path groupoid is isomorphic to $G_T$, as follows:
	\begin{itemize}
		\item $\operatorname{Obj}(\Lambda_T) \coloneqq \Omega$;
		\item $\operatorname{Mor}(\Lambda_T)\coloneqq \{(x,n,y) \in \Omega\times\N^k\times \Omega : T^n x = y\}$;
		\item $d(x,n,y) = n$, $r(x,n,y) = x$ and $s(x,n,y) = y$.
	\end{itemize}
	It is not hard to check that this satisfies factorisation and gives a source-free proper topological $2$-graph. The boundary path groupoid $G_{\Lambda_T}$ in the sense of Yeend \cite{Yeendtopk} has as units the space of infinite paths $\Lambda_T^\infty$ of the topological $2$-graph since it is source-free and proper (see \cite[p.~1448]{ArmBrown2018}). Furthermore, the infinite-path space is homeomorphic to $\Omega$ via the map $\Lambda_T^\infty \to \Omega$ that sends an infinite path $\zeta$ to its range vertex $\zeta(0)\in \Omega$.  After identifying $\Omega$ and  $\Lambda_T^\infty$,  the boundary-path groupoid $G_{\Lambda_T}$  equals the  Deaconu--Renault groupoid $G_T$.
	
We use the construction in \cite[\S5]{CLSV} to build a product system:
	\begin{itemize}
		\item For each $m \in \N^2$, let $\Lambda_T^m \coloneqq d^{-1}(m) = \{(x,m,y) \in \Lambda_T : T^mx = y\}$;
		\item For each $m\in \N^2$, define $E_{T^m}$ to be the right-Hilbert $C_0(\Omega)$-bimodule associated to the surjective local homeomorphism $T^m \colon  \Omega\to \Omega$ as defined in \S\ref{sec:KK}. Note that $E_{T^m}$ is the topological-graph module associated to the topological ($1$-graph) $(\Omega,\Omega,\operatorname{id},T^m)$ which is isomorphic to the topological graph $(\Lambda^0,\Lambda^m,r|_{\Lambda^m},s|_{\Lambda^m})$. Hence, writing $\mathbf{X}_m \coloneqq E_{T^m}$, we can apply all of \cite[Sec.~5]{CLSV} to the family $\mathbf{X_T}\coloneqq \bigsqcup_{m\in\N^2}\mathbf{X}_m$;
		\item For $f\in C_c(\Omega)\subset \mathbf{X}_m$ and $g\in C_c(\Omega)\subset \mathbf{X}_n$, define $fg \colon  \Omega \to \C$ by $(fg)(x) = f(x)g(T^mx)$; this $fg$ lives in $C_c(\Omega) \subseteq \mathbf{X}_{n+m}$; we then extend this formula by continuity to a multiplication $(\xi,\eta) \mapsto \xi\eta$ from $X_m \times X_n$ to $X_{m+n}$.
	\end{itemize}
	Under this multiplication, the family $\mathbf{X_T}\coloneqq \bigsqcup_{n\in\N^k}X_n$ of right-Hilbert $C_0(\Omega)$-bimodules is a compactly aligned product system over $\N^2$, and according to \cite[Theorem~5.20]{CLSV}, the Cuntz--Nica--Pimsner algebra $\NO_{\mathbf{X_T}}$ associated to $\mathbf{X_T}$ is isomorphic to the $C^*$-algebra $C^*(G_{\Lambda_T})$ of the boundary-path groupoid and hence isomorphic to $C^*(G_T)$. The isomorphism
	\[S'\colon  \NO_{\mathbf{X_T}} \to  C^*(G_T)\] works as follows: for each $m\in \N^2$, there exists a map $\psi_m \colon  \mathbf{X}_m \to C^*(G_T)$ such that, for $f \in C_c(\Omega) \subset \mathbf{X}_m$ and $(x,p,y)\in G_T$,
		\[\psi_m(f)(x,p,y) = \begin{cases}
		f(x) \text{ if } p = m \text{ and } T^mx = y, \\
		0 \text{ otherwise};
		\end{cases}\]
		and then $S'$ maps $j_\mathbf{X_T}(f)$ to $\psi_m(f)$ (see \cite[Theorem 5.20]{CLSV}: our map $\psi_m$ is the composition of the map from the top left to the bottom right in the diagram of \cite[Theorem 5.20]{CLSV} with the canonical inclusion of $\mathbf{X}_m$ in $\TT_{\operatorname{cov}}(\mathbf{X_T})$).
	
\subsection*{Applying the Deaconu--Fletcher construction}
We now apply the Deaconu--Fletcher construction discussed in \S\ref{sec:fletcher} to the product system above, with
\[
A \coloneqq C_0(\Omega), X_1 \coloneqq \mathbf{X}_{(1,0)}=E_{T_1} \text{ and \ }X_2 \coloneqq \mathbf{X}_{(0,1)}=E_{T_2}.
\]
Then $Y_1 \coloneqq X_2\otimes_A\OO_{X_1}$ is a right-Hilbert $\OO_{X_1}$-bimodule (and not only an $A$--$\OO_{X_1}$-bimodule), and $\OO_{Y_1} \cong \NO_\mathbf{X_T} \cong C^*(G_T)$ by Theorem~\ref{thm:flecther}. Furthermore, the composition
	\begin{equation}\label{eq:isomorphismS}
	S\colon \OO_{Y_1}\to  C^*(G_T)
	\end{equation}
	of these isomorphisms satisfies the following: for $a\in A$, $\xi \in C_c(\Omega)\subset X_1$ and $\eta \in C_c(\Omega)\subset  X_2$, define $\eta\xi \colon \Omega \to \C$ by $(\eta\xi)(x) = \eta(x)\xi(Tx)$; then by Remark \ref{rmk:expfletcheriso},
		\begin{align}\label{eq:formulasforS}
		S(j_1(j_1(a))) & = \psi_{0}(a), \notag\\
		S(j_1(j_2(\xi))) & = \psi_{(1,0)}(\xi),\notag \\
		S(j_2(\eta \otimes j_1(a))) & = \psi_{(1,0)}(\eta \cdot a),\text{ and}\notag \\
		S(j_2(\eta\otimes j_2(\xi))) & =	\psi_{(1,1)}(\eta \xi).	
		\end{align}
		The first of these equations will be particularly important later, so we expand on our shortened notation for future reference: the map $\psi_0 \colon  A \to C^*(G_T)$ is given by
		\[
		\psi_0(a)(x,p,y) = \begin{cases}
		a(x) &\text{ if } p = 0 \text{ and } x=y, \\
		0 &\text{ otherwise},
		\end{cases}
		\]
		and then
		\begin{equation}\label{eq:phisub0}
		S\big(j_1^{(2)}\big(j_1^{(1)}(a)\big)\big) = \psi_0(a) \in C^*(G_T).
		\end{equation}
	
\subsection*{Invariant ideals and subsets}
Let $H$ be an open subset of $\Omega$ that is invariant for the dynamics	
and consider the ideal $I\coloneqq C_0(H)$ of $A=C_0(\Omega)$.   The two possible definitions of ``invariant'' coincide:

\goodbreak		
	
\begin{lemma}\label{lem:invarianceequivalence}
Let $H\subset \Omega$. Then $H$ is  invariant under pre-images and images of $T_1$ and $T_2$ if and only if $H$ is  an invariant subset of the unit space $\Omega$ of $G_T$, that is, $r(s^{-1}(H))=H$.
\end{lemma}

\begin{proof}
First suppose that $H$ is  invariant under pre-images and images of $T_1$ and $T_2$.  We always have $H\subset r(s^{-1}(H))$.  Let $y\in r(s^{-1}(H))$. Then there exists $(y, p-q, x)\in G_T$ with $x\in H$. Then $T^py=T^qx$ implies $y\in (T^p)^{-1}(T^q(H))\subset H$. Thus $r(s^{-1}(H))=H$.

Suppose that $r(s^{-1}(H))=H$. Let $x\in H$. Then $(T_ix, -e_i, x)\in G_T$ implies $T_ix\in r(s^{-1}(H))=H$, giving $T_i(H)\subset H$. Let $y\in T_i^{-1}(H)$, say $T_i(y)=x\in H$. Then $(y, e_i, x)\in G_T$ implies $y\in r(s^{-1}(H))=H$, giving $T_i^{-1}(H)\subset H$.
\end{proof}

	\begin{lemma}\label{lem:invariantideal}
		Let $H$ be an open invariant subset of the unit space $\Omega$ of $G_T$ and let $i\in\{1,2\}$.
		\begin{enumerate}
			\item\label{lem:invideal1} $C_0(H)$ is an $X_i$-invariant ideal of $C_0(\Omega)$; that is $C_0(H)X_i\subset X_iC_0(H)$.
			\item\label{lem:invideal2}  The inclusion $\iota_H \colon  C_c(H) \to C_c(\Omega)$ extends to an isomorphism of the right-Hilbert $C_0(H)$-bimodule $X_i^H\coloneqq E_{T_i|_H}$ onto $X_i  C_0(H)$.
		\end{enumerate}
	\end{lemma}

	\begin{proof} By Lemma~\ref{lem:invarianceequivalence},  $H$ is  invariant under pre-images and images of $T_1$ and $T_2$.
		Fix $f\in C_0(H)$ and $\xi \in X_i$. We must show that $f \cdot \xi \in X_i  C_0(H)$. By linearity and continuity, it suffices to consider $f = 1_U$ for some $U\subset H$ for which $T_i|_U \colon  U\to T_i(U)$ is a homeomorphism of compact open sets. Then $1_{T_i(U)} \in C_0(H)$ because $H$ is invariant under $T_i$, and so
		\[
			1_U \cdot \xi = (1_U \cdot \xi)\cdot 1_{T_i(U)} \in X_i  C_0(H).
		\]
		Thus $C_0(H)$ is $X_i$-invariant, proving~(\ref{lem:invideal1}).

        Write $\operatorname{CO}_T(\Omega)$ (respectively $\operatorname{CO}_T(H)$) for the collections of compact open subsets $U$ of $\Omega$ (respectively $H$) such that $T_1|_U$ and $ T_2|_U$ are homeomorphisms onto their ranges. Fix $U, V \in \operatorname{CO}_T(H)$. Then
		\[
		\langle \iota_H(1_U), \iota_H(1_V)\rangle_{C_0(\Omega)}
			= 1_{T_i(U) \cap T_i(V)}
			= \iota_H(\langle 1_U, 1_V\rangle_{C_0(H)}).
		\]
		So by linearity and continuity, $\iota_H$ extends to an isometric linear map, also denoted $\iota_H$, from $X_i^H$ to $X_i$.
		For $U \in \operatorname{CO}_T(H)$, we have $1_U = 1_U \cdot 1_{T_i(U)} \in X_i C_0(H)$, and so the range of $\iota_H$ is contained in $X_i C_0(H)$. Furthermore, for $U \in \operatorname{CO}_T(\Omega)$ and $V \in \operatorname{CO}_T(H)$, we have $1_U \cdot \iota_H(1_V) = 1_{U\cap T^{-1}(V)}$. Hence $C_c(\Omega)  \iota_H(C_c(H)) = \iota_H(C_c(H))$ and since $C_c(\Omega)  \iota_H(C_c(H))$ is dense in $X_i  C_0(H)$, we deduce that $\iota_H(X_i^H) = X_i C_0(H)$. Checking that $\iota_H$ respects the left and right actions is straightforward, proving~(\ref{lem:invideal2}).	
\end{proof}

Let $H\subset\Omega$ be an open invariant subset of the unit space of $G_T$.  Then the restriction
\[
(G_T)|_H=\{\gamma\in G_T: r(\gamma), s(\gamma)\in H\}
\]
is a locally compact Hausdorff ample groupoid. It follows from Lemma~\ref{lem:invarianceequivalence} that the groupoid $(G_T)|_H$ coincides with the Deaconu--Renault groupooid $G_{T|_H}$ of the restricted commuting maps $T_1|_H$ and $T_2|_H$.
Further, let $i\colon  C^*((G_T)|_H) \to C^*(G_T)$  be the homomorphism induced by  inclusion  and extension by $0$ of $C_c((G_T)|_H)$ in $C_c(G_T)$. Then $C^*((G_T)|_H)$ is isomorphic to an ideal of $C^*(G)$ and there is an exact sequence
	\begin{equation}\label{eq:exactsequencegroupoids}
		\begin{tikzcd}
			0 \arrow[r] & C^*((G_T)|_H) \arrow[r,"i"] & C^*(G_T) \arrow[r] & C^*((G_T)|_{\Omega\setminus H}) \arrow[r] & 0
		\end{tikzcd}
\end{equation}
by \cite[Theorem~5.1]{tfb^3}.

In Lemma~\ref{lem:invariantideal} we set $X_i^H\coloneqq E_{T_i|_H}$.  Now we also write
$Y_1^H$ for the right-Hilbert $C_0(H)$-bimodule $(X_2 C_0(H))\otimes_{C_0(H)}\OO_{X_1 C_0(H)}$ from the Deaconu--Fletcher construction of  \S\ref{sec:fletcher} and  $Y(H)_1$ for the right-Hilbert $C_0(H)$-bimodule $X_2^H \otimes_{C_0(H)} \OO_{X_1^H}$. By Lemma~\ref{lem:invariantideal}, $Y_1^H$ and $Y(H)_1$ are isomorphic, and hence there is an isomorphism $\OO_{Y_1^H} \cong \OO_{Y(H)_1}$; we identify the two $C^*$-algebras via this isomorphism. Furthermore, by Lemmas~\ref{lem:invariantideal}(\ref{lem:invideal1}) and \ref{lem:fletchermorphisms}, there exists a canonical homomorphism $\phi\colon\OO_{Y_1^H}\to \OO_{Y_1}$. We now check that $\phi$ corresponds to the inclusion $i$ of $C^*((G_T)|_H)$  in $C^*(G_T)$.

	\begin{lemma}\label{lem:phiandi} Let $H$ be an open invariant subset of the unit space $\Omega$ of $G_T$ and
		let $i\colon C^*((G_T)|_H) \to C^*(G_T)$ be the homomorphism of \eqref{eq:exactsequencegroupoids}. The construction of the isomorphism $S\colon \OO_{Y_1}\to C^*(G_T)$ above applied to the reduced system $T|_H$ gives an isomorphism $S_H \colon  \OO_{Y(H)_1} \to C^*((G_T)|_H)$,
		and
		the diagram
		\[
		\begin{tikzcd}
			\OO_{Y_1^H} \arrow[d, "\phi"] \arrow[r, "\cong"] & \OO_{Y(H)_1} \arrow[r, "S_H"] & C^*((G_T)|_H) \arrow[d, "i"] \\
			\OO_{Y_1} \arrow[rr, "S"]                      & & C^*(G_T)
		\end{tikzcd}\]
		commutes.
	\end{lemma}

	\begin{proof}
		The Deaconu--Renault groupoid $(G_T)|_H$ gives rise to $Y(H)_1$ in the same way that $G_T$ gives rise to $Y_1$, and so we obtain $S_H \colon  \OO_{Y(H)_1}\to C^*((G_T)|_H)$ in the same way as before. 
		Fix $\eta\in C_c(H) \subset X_2^H,\xi \in C_c(H)\subset X_1^H$. Temporarily write $\overline{\eta}$ and $\overline{\xi}$ for the extensions of $\eta$ and $\xi$ to elements of $C_c(\Omega)$ that vanish on $\Omega \setminus H$, regarded as elements of $X_2$ and $X_1$ (that is, $\iota_H(\eta)$ and $\iota_H(\xi)$ as in Lemma \ref{lem:invariantideal}(\ref{lem:invideal2}) respectively).
		Then
		\[
		S(\phi(j_2(\eta\otimes j_2(\xi)))) = S(j_2(\overline{\eta}\otimes j_2(\overline{\xi}))) = \psi_{(1,1)}(\overline{\eta}\,\overline{\xi}).
		\]
		 The product $\overline{\eta}\,\overline{\xi} \in \mathbf{X}_{(1,1)}$ in the product system $\mathbf{X_T}$ is given by $(\overline{\eta}\,\overline{\xi})(x) = \overline{\eta}(x)\overline{\xi}(T_1 x)$, so has support in $H$. In particular, $\overline{\eta}\,\overline{\xi} = \overline{\eta\xi}$. Hence (see the formulas at \eqref{eq:formulasforS})
		 \[
		 S(j_2(\overline{\eta}\otimes j_2(\overline{\xi}))) = \psi_{(1,1)}(\overline{\eta\xi}) = i(\psi^H_{(1,1)}(\eta\xi)).
		 \]
		 Since $\psi^H_{(1,1)}(\eta\xi) = S_H(j_2(\eta\otimes j_2(\xi)))$ by definition, we deduce that \[i(S_H(j_2(\eta\otimes j_2(\xi)))) = S(\phi(j_2(\eta\otimes j_2(\xi)))).\] Similar calculations show that $i \circ S_H$ and $S \circ \phi$ agree on elements of the form $j_2(\eta\otimes j_1(a))$, $j_1(j_2(\xi))$ and $j_1(j_1(a))$. Since these elements generate $\OO_{Y(H)_1}$, the result follows.
	\end{proof}

\subsection*{The main theorem of the section}
 We have now set up the background required  to pull over Theorem~\ref{thm:naturalses} to obtain information about the K$_0$ groups of $G_T$ and $G_T|_H$. Theorem~\ref{thm:DRkgroups2} below is the analogue of \cite[Theorem~3.20]{aHNS} for the $C^*$-algebra of  a $2$-graph and  a gauge-invariant ideal.  While the theorems are similar, the approaches taken to prove them  have been very different, as discussed in the introduction. With Theorem~\ref{thm:DRkgroups2} in hand we will be able to follow the programme in \cite{aHNS} to obtain a result about stable finiteness of  $C^*(G_T)$  in \S\ref{sec:SFE} below.

Consider $[X_i]\colon K_0(C_0(\Omega))\to K_0(C_0(\Omega))$,
 and let  $\sigma \colon C_c(\Omega, \Z) \to K_0(C_0(\Omega))$ be the isomorphism characterised by  $\sigma(1_U) = [1_U]_0$ for compact open $U \subset \Omega$.  Then Theorem~\ref{thm:tdlcKtheory} with $X=Y=\Omega$ and $T = T_i$ implies  that
\begin{equation*}
1-[X_i]\circ \sigma = \sigma \circ (1 - (T_i)_*).
\end{equation*}
It follows that
$\sigma$ descends to an isomorphism
	\begin{equation}\label{eq:kers}
	\tilde{\sigma} \colon  \coker(1-(T_1)_*,1-(T_2)_*) \to \coker(1-[X_1],1-[X_2])
	\end{equation}
	and restricts to an isomorphism
	\begin{equation}\label{eq:cokers}
	\sigma| \colon  \ker\begin{pmatrix} 1-(T_1)_* \\ 1- (T_2)_* \end{pmatrix} \to \ker{\begin{pmatrix}1-[X_1]\\ 1-[X_2]\end{pmatrix}}.
	\end{equation}

	\begin{theorem}\label{thm:DRkgroups2}
		Let $T_1,T_2$ be commuting surjective local homeomorphisms on a second-countable totally disconnected locally compact Hausdorff space $\Omega$. Let $G_T$ denote the associated rank-$2$ Deaconu--Renault groupoid  and let $H$ be an invariant subset of its unit space.  Let $\imath\colon C_0(H)\to C_0(\Omega)$ be the inclusion, let $\sigma\colon C_c(\Omega,\Z)\to K_0(C_0(\Omega))$ be the homomorphism from \eqref{eq:sigma def} and let $i\colon C^*((G_T)|_H)\to C^*(G_T)$ be the inclusion from \eqref{eq:exactsequencegroupoids}.
		There are homomorphisms $j_\Omega,j_H, \tau_\Omega, \tau_H$ such that
		 the diagram
		\[
		\begin{tikzcd}
			0 \arrow[r] &[-1em] {\coker(1-(T_1|_H)_*,1-(T_2|_H)_*)} \arrow[r, "j_H"]
			\arrow[d, "\tilde{\sigma}^{-1}\circ\tilde\imath\circ \tilde{\sigma}_H"]
			&[-1ex] K_0(C^*((G_T)|_H)) \arrow[r, "\tau_H"] \arrow[d, "i_*"] &[-1em] \ker\begin{pmatrix}1-(T_1|_H)_* \\ 1- (T_2|_H)_* \end{pmatrix} \arrow[r] \arrow[d, "\sigma|^{-1}\circ\imath_*\circ{\sigma_H}|"] & 0 \\
			0 \arrow[r] & {\coker(1-(T_1)_*,1-(T_2)_*)} \arrow[r, "j_\Omega"]                                     & K_0(C^*(G_T)) \arrow[r, "\tau_\Omega"]                            & \ker\begin{pmatrix} 1-(T_1)_* \\ 1- (T_2)_* \end{pmatrix} \arrow[r]                          & 0
		\end{tikzcd}
		\]
		commutes and has exact rows. Moreover, for any compact open set $U \subset \Omega$,
		\begin{align}
		&j_\Omega(1_U + \image(1-(T_1)_*,1-(T_2)_*)) = [\psi_0(1_U)]_0\label{eq:j_Omega}\quad \text{(and likewise for $j_H$); and}
		\\
		&\sigma|^{-1}\circ\imath_*\circ \sigma|_H(1_U + \image(1-(T_1|_H)_*,1-(T_2|_H)_*)=\imath(1_U)+ \image(1-(T_1)_*,1-(T_2)_*).\label{eq:imathpulledover}
		\end{align}
Finally, $\sigma|^{-1}\circ\imath_*\circ \sigma_H|$ is injective.
	\end{theorem}
	\begin{proof}
Consider  the diagram below.  The inner two exact rows connected by the vertical maps are obtained by  applying Theorem~\ref{thm:naturalses} with  $I \coloneqq C_0(H) \lhd A\coloneqq C_0(\Omega)$.
Let $\tilde\sigma$ and $\sigma|$  be the isomorphisms from Equations \eqref{eq:kers}~and~\eqref{eq:cokers},  and let  $S_*$ be the isomorphism  induced from the isomorphism $S$ at \eqref{eq:isomorphismS}. Then  we augment the inner two rows using the homomorphisms of the previous sentence as the vertical maps as shown:
	\[
		\begin{tikzcd}
		0 \arrow[r] &[-1em] {\coker(1-(T_1|_H)_*,1-(T_2|_H)_*)} \arrow[r, "j_H"] \arrow[d, "\tilde\sigma_H"] & K_0(C^*((G_T)|_H)) \arrow[r, "\tau_H"] & \ker\begin{pmatrix}1-(T_1|_H)_* \\ 1- (T_2|_H)_* \end{pmatrix} \arrow[r] \arrow[d, "\sigma_H|"]  &[-1em] 0 \\
		0 \arrow[r] & {\coker(1-[X_1\cdot I],1-[X_2\cdot I])} \arrow[r, "j_I"] \arrow[d, "\tilde{\imath}"] & K_0(\OO_{Y_1^I}) \arrow[u, "(S_H)_*"] \arrow[r, "\tau_I"] \arrow[d, "\phi_*"] & {\ker{\begin{pmatrix}1-[X_1\cdot I]\\ 1-[X_2\cdot I]\end{pmatrix}}} \arrow[r] \arrow[d, "\imath_*"] & 0 \\
		0 \arrow[r] & {\coker(1-[X_1],1-[X_2])} \arrow[r, "j"] & K_0(\OO_{Y_1}) \arrow[d, "S_*"] \arrow[r, "\tau"] & {\ker{\begin{pmatrix}1-[X_1]\\ 1-[X_2]\end{pmatrix}}} \arrow[r] & 0 \\
	    0 \arrow[r] & {\coker(1-(T_1)_*,1-(T_2)_*)} \arrow[r, "j_\Omega"] \arrow[u, "\tilde\sigma"] & K_0(C^*(G_T)) \arrow[r, "\tau_\Omega"] & \ker\begin{pmatrix} 1-(T_1)_* \\ 1- (T_2)_* \end{pmatrix} \arrow[r] \arrow[u, "\sigma|"] & 0.\\
		\end{tikzcd}
	\]
	We then define $j_\Omega, j_H, \tau_\Omega$  and $\tau_H$ to be the unique homomorphisms that make the whole diagram commute. Then the bottom and top rows are exact by construction.  Notice  that $i=S_*\circ\phi_*\circ (S_H^{-1})_*$ by Lemma~\ref{lem:phiandi}.
	
Fix a compact open set $U \subset \Omega$.  Then
	\begin{align*}
		j_\Omega(1_U + \image(1-(T_1)_*,1-(T_2)_*))
			&= S_* \circ j \circ \tilde{\sigma}\big(1_U + \image(1-(T_1)_*,1-(T_2)_*)\big)\\
			&= S_* \circ j \big(1_U + \image(1-[X_1],1-[X_2])\big)\\
			&= S_* \big(\big[(j_1^{(2)} \circ j_1^{(1)})_*(1_U)\big]_0\big)\\
			&= \big[S\big(j_1^{(2)} \circ j_1^{(1)}(1_U)\big)\big]_0\\
			&= \big[\psi_0(1_U)\big]_0
	\end{align*}
	by~\eqref{eq:phisub0}. The computation for $j_H$ is identical.
Similarly,
\begin{align*}
\sigma|^{-1}\circ\imath_*\circ \sigma_H|(1_U + \image(1-(T_1|_H)_*,1-(T_2|_H)_*)
&=\sigma|^{-1}(\imath_*(1_U)+ \image(1-[X_1I],1-[X_2I])\\
&=\imath(1_U) +  \image(1-(T_1)_*,1-(T_2)_*)
\end{align*}
Finally, since $K_0(I) = C_c(H,\Z)\hookrightarrow C_c(\Omega,\Z) = K_0(A)$, the homomorphism $\iota_*\colon K_0(I)\to K_0(A)$ is injective, and so the restriction
		\[
		\imath_*\colon \ker\begin{pmatrix}1-[X_1\cdot I]\\1-[X_2\cdot I]\end{pmatrix}\to \ker\begin{pmatrix}1-[X_1]\\1-[X_2]\end{pmatrix}
		\]
appearing in Theorem~\ref{thm:naturalses} is injective. Thus $\sigma|^{-1}\circ\imath_*\circ \sigma_H|$ is injective as well.
	\end{proof}

\section{Application: stable finiteness of extensions of Deaconu--Renault groupoid \texorpdfstring{$C^*$}{C*}-algebras}\label{sec:SFE}
	
	Let $G$ be an Hausdorff, ample groupoid and let $\mathcal{C}$ denote the family of all compact open bisections in $G$. As defined in \cite[Definition 6.4]{RainoneSims2020}, the \textit{coboundary subgroup} of $G$ is the subgroup
	$$H_G \coloneqq \big\langle 1_{s(E)} - 1_{r(E)} : E\in \mathcal{C}\big\rangle$$
	of $C_c(G^{(0)},\Z)$. The \textit{coboundary condition} is then defined
	to be
	\begin{align}\tag{C}\label{C}
        H_G \cap C_c(G^{(0)},\N) = \{0\}.
    \end{align}

	Corollary~6.6 of \cite{RainoneSims2020} shows that if $G$ is a minimal ample groupoid, then the coboundary condition characterises stable finiteness of $C^*_r(G)$  (see also \cite[Theorem~5.14]{BonickeLi} for a similar result). As an application of our $K$-theory calculations, we extend this to the situation where a Deaconu--Renault groupoid associated to an action of $\N^2$ by local homeomorphism has a unique nontrivial open invariant subset (thus is not minimal).

	Corollary~6.6 of \cite{RainoneSims2020} can be thought of as a groupoid analogue of \cite[Theorem~1.1]{CaHS} which states that for a row-finite cofinal $k$-graph with no sources, stable-finiteness of the associated $C^*$-algebra is equivalent to the adjacency matrices of the graph satisfying a so-called matrix condition. Cofinality of a $k$-graph corresponds to minimality of its groupoid, and so we might expect that for a Deaconu--Renault groupoid,  the coboundary condition \eqref{C} is related to a matrix condition analogous to the one in \cite{CaHS}; we prove this below.

Let $T_1,T_2$ be commuting surjective local homeomorphisms on a second-countable totally disconnected locally compact Hausdorff space $\Omega$, and let $G_T$ be the associated rank-$2$ Deaconu--Renault groupoid defined at \eqref{eq-DRgroupoid}.
For the triple $(\Omega,T_1,T_2)$, define the \textit{generalised matrix condition} as
	\begin{align}\label{M}\tag{M}
		\{(1-(T_1)_*)f + (1-(T_2)_*)g : f,g \in C_c(\Omega,\Z)\}\cap C_c(\Omega,\N) = \{0\}.
	\end{align}
In Proposition~\ref{prop: matrixconditioncoboundary} below we prove that $G_T$ satisfies the coboundary condition  if and only if $(\Omega,T_1,T_2)$ satisfies the generalised matrix condition.  We start with the following observation.
	
	\begin{lemma}\label{lem:uniondisjointbisections}
	There are compact open subsets  $U_i, V_i\subset \Omega$ and elements $p_i, q_i\in\N^2$  (indexed by $i \in \mathbb{N}$) such that $T^{p_i}|_{U_i}$, $T^{q_i}|_{V_i}$ are injective and $T^{p_i}({U_i})=T^{q_i}({V_i})$ and such that
	\[G_T=\bigsqcup_{i\in\N} Z(U_i, p_i, q_i, V_i).\]
	\end{lemma}
	
	\begin{proof}
	The collection of sets $Z(U, p, q, V)= \{(x,p-q,y):\text{$x\in U$, $y\in V$ and $T^{p}x=T^{q}y$}\}$,  where $U,V$ are compact open subsets of $\Omega$, $p,q\in\N^2$, and  $T^{p_i}|_{U_i}$ and $T^{q_i}|_{V_i}$ are injective and $T^{p_i}({U_i})=T^{q_i}({V_i})$, constitute a basis of compact open sets for a second-countable  topology of $G_T$ by \cite[Lemma~3.1]{SimsWilliams2016}. Thus we can cover $G_T$ with a countable collection $\{Z_i\colon i\in\N\}$ of such sets. For $n \ge 0$, let $E_n := Z_n\setminus\bigcup_{0 \le i < n} Z_i$ to obtain a disjoint cover. Since each $E_n$ is the relative complement of one compact open set in another, it is itself compact open. Since $T^{p_n}|_{U_n}$ and $T^{q_n}|_{V_n}$ are homeomorphisms onto a common compact open set $W_n$, the map $z \mapsto \big((T^{p_n}|_{W_n})^{-1}(z), p_n - q_n, (T^{q_n}|_{W_n})^{-1}(z)\big)$ is a homeomorphism of $W_n$ onto $Z_n$. The pre-image $W'_n$ of $E_n$ under this homeomorphism is compact open, and $U'_n := (T^{p_n}|_{W_n})^{-1}(W'_n)$ and $V'_n := (T^{q_n}|_{W_n})^{-1}(W'_n)$ are compact open sets such that $T^{p_n}|_{U'_n}$ and $T^{q_n}|_{V'_n}$ are injective and $T^{p_n}({U'_n}) = W'_n = T^{q_n}({V'_n})$. Since $E_n = Z(U'_n, p_n, q_n, V'_n)$, this proves the lemma.
	\end{proof}
	
\begin{proposition}\label{prop: matrixconditioncoboundary}
Let $T_1,T_2$ be commuting surjective local homeomorphisms of a second-countable totally disconnected locally compact Hausdorff space $\Omega$. Let $G_T$ be the associated rank-$2$ Deaconu--Renault groupoid. Then
\begin{equation}\label{eq:coboundaryandmatrixcondition}
		\{(1-(T_1)_*)f + (1-(T_2)_*)g : f,g \in C_c(\Omega,\Z)\} = \big\langle 1_{s(E)} - 1_{r(E)} : E\in \mathcal{C}\big\rangle,
		\end{equation}
and thus $G_T$ satisfies the coboundary condition \eqref{C} if and only if $(\Omega,T_1,T_2)$ satisfies the generalised matrix condition \eqref{M}.
\end{proposition}
\begin{proof}
	Fix $i\in\{1,2\}$ and fix a compact neighbourhood $U\subset\Omega$ such that $T_i|_{U}$ is a homeomorphism.  Let $E=Z(T_i(U),0,e_i,U)$. Then $E\in \mathcal{C}$ because $T_i$ is injective on $U$.  We have
 \[(1-(T_i)_*)(1_U) = 1_U - 1_{T_i(U)} = 1_{s(E)} - 1_{r(E)}.\]
 Since  $C_c(\Omega,\Z)$ is generated by characteristic functions of compact open sets $U$ on which $T_i$ is a homeomorphism, this shows that the left-hand side of \eqref{eq:coboundaryandmatrixcondition} is contained in the right-hand side.

For the other inclusion, fix  $E\in\mathcal{C}$.  By Lemma~\ref{lem:uniondisjointbisections} we can write $E$ as a disjoint finite union $E=\bigsqcup_i Z(U_i, p_i, q_i, V_i)$ where all the $U_i, V_i$ are compact open,  $T^{p_i}|_{U_i}$, $T^{q_i}|_{V_i}$ are injective and $T^{p_i}({U_i})=T^{q_i}({V_i})$. Hence
\[1_{s(E)}-1_{r(E)} = \sum_i 1_{U_i} - 1_{V_i}.\]
Since the left-hand side of \eqref{eq:coboundaryandmatrixcondition}  is a group, it suffices to show that it contains each $1_{U_i} - 1_{V_i}$. So fix $E=Z(U, p, q, V)$  as above. Since $T^p(U)=T^q(V)$ we have $1_{U} - 1_{V}=1_{U} -1_{T^p(U)}+1_{T^q(V)}- 1_{V}$. Thus it suffices to show that $1_{U} -1_{T^p(U)}$ is in the left-hand side of \eqref{eq:coboundaryandmatrixcondition}.  We show this by induction on $|p|=p_1+p_2$.  If $|p|=1$, then $p=e_i$ for $i\in\{1,2\}$ and  $1_{U} -1_{T^p(U)}=(1-(T_i)_*) (1_U)$, as needed. Now assume that if $|p|\leq N$, then $1_{U_i} -1_{T^p(U)}$ is in the left-hand side \eqref{eq:coboundaryandmatrixcondition}.  Fix $p$ with $|p|=N+1$.  Write $p=q+e_i$ so that $|q|=N$. Then
\[
1_{U} -1_{T^p(U)}=1_{U} -1_{T^q(U)}+1_{T^q(U)}-1_{T_i(T^q(U))}
\]
which is in the left-hand side \eqref{eq:coboundaryandmatrixcondition}  by the induction hypothesis and the base case.  Thus \eqref{eq:coboundaryandmatrixcondition}  holds.
\end{proof}

\begin{remark} The proof of Proposition \ref{prop: matrixconditioncoboundary} generalises to  $n$ commuting local homeomorphisms.	\end{remark}
	
We have the following partial analogue of \cite[Lemma~4.1]{aHNS}.
	
	\begin{lemma}\label{lem:MforH}
Let $H\subset \Omega$ be an open invariant subset of the unit space of $G_T$.	
		Let \[\tilde\sigma|^{-1}\circ \tilde{\imath}\circ\tilde\sigma_H \colon  \coker(1-(T_1|_H)_*,1-(T_2|_H)_*) \to \coker(1-(T_1)_*,1-(T_2)_*)\]  be the homomorphism from \eqref{eq:imathpulledover}.
		If $(\Omega,T_1,T_2)$ satisfies \eqref{M}, then $(H,T_1|_H,T_2|_H)$ satisfies \eqref{M} and
		\begin{equation*}
			\ker(\tilde\sigma^{-1}\circ \tilde{\imath}\circ\tilde\sigma_H)\cap\big [C_c(H,\N) +  \image\big(1-(T_1|_H)_*,1-(T_2|_H)_*\big)\big] = \{0\}.
		\end{equation*}
	\end{lemma}

	\begin{proof}
Suppose that $(\Omega,T_1,T_2)$ satisfies \eqref{M}. Then
\begin{align*}
			\{(1-(T_1|_H)_*)f + & (1-(T_2|_H)_*)g : f,g \in C_c(H,\Z)\}\cap C_c(H,\N) \\
			& \subset \{(1-(T_1)_*)f + (1-(T_2)_*)g : f,g \in C_c(\Omega,\Z)\}\cap C_c(\Omega,\N) = \{0\}.
\end{align*}
Hence $(H,T_1|_H,T_2|_H)$ satisfies \eqref{M}. Furthermore,
\begin{align*}
\ker(\tilde\sigma^{-1}\circ \tilde{\imath}\circ\tilde\sigma_H )\cap {}& \big(C_c(H,\N) + \image(1-(T_1|_H)_*,1-(T_2|_H)_*)\big) \\
& \subset \{(1-(T_1)_*)f + (1-(T_2)_*)g : f,g \in C_c(\Omega,\Z)\}\cap C_c(\Omega,\N) = \{0\}.\qedhere
\end{align*}
\end{proof}
	
	We can now state and prove our main theorem on stably finite extensions. It is a Deaconu--Renault groupoid analogue of \cite[Theorem~4.6]{aHNS} and follows its proof closely.
	
	\begin{theorem}\label{thm:lastmain}
		Let $T_1,T_2$ be commuting surjective local homeomorphisms on a  second-countable totally disconnected locally compact Hausdorff space $\Omega$.  Let $G=G_T$ denote the associated rank-$2$ Deaconu--Renault groupoid  and let $H\subset\Omega$ be an open invariant subset of its unit space.
		Let \[\kappa_H \colon  C_c(H,\Z) \to K_0(C^*(G|_H))\] be the composition of the quotient map $q\colon C_c(H,\Z) \to \coker(1-(T_1|_H)_*,1-(T_2|_H)_*)$ and the homomorphism $j_H \colon  \coker(1-(T_1|_H)_*,1-(T_2|_H)_*) \to K_0(C^*(G|H))$ from \eqref{eq:j_Omega}.
		Assume that
		\begin{align*}
			\{(1-(T_1)_*)f + (1-(T_2)_*)g : f,g \in C_c(X,\Z)\}\cap C_c(X,\N) = \{0\}, \tag{M}\text{ and}
		\end{align*}
		\begin{align}\label{P}\tag{P}
		\kappa_H(C_c(H,\Z)) \cap K_0(C^*(G|_H))_+ = \kappa_H(C_c(H,\N)).
		\end{align}
	If $C^*(G|_H)$ and $C^*(G|_{\Omega\setminus H})$ are stably finite, then $C^*(G)$ is stably finite.
	\end{theorem}
	
\begin{proof}
Since $H$ is an invariant subset of the unit space of $G_T$,  \cite[Theorem~5.1]{tfb^3} gives a short exact sequence
		\[
		\begin{tikzcd}
			0 \arrow[r] & C^*(G|_H) \arrow[r,"i"] & C^*(G) \arrow[r] & C^*(G|_{\Omega\setminus H}) \arrow[r] & 0.
		\end{tikzcd}
		\]
 Since $C^*(G|_H)$ and $C^*(G|_{\Omega\setminus H})$ are stably finite by assumption, by \cite[Lemma 1.5]{Spielberg1988}  $C^*(G)$ is stably finite if and only if 		
        \begin{align}\label{S}\tag{S}
			\ker(i_*)\cap K_0(C^*(G|_H))_+ = \{0\}.
		\end{align}	
		
Using assumption \eqref{M}, by Lemma \ref{lem:MforH}  we have
\begin{equation*}
		\ker(\tilde\sigma^{-1}\circ \tilde{\imath}\circ\tilde\sigma_H)\cap \big(C_c(H,\N) + \image(1-(T_1|_H)_*,1-(T_2|_H)_*)\big) = \{0\}.
\end{equation*} 	
Consider the left square of the commuting diagram in Theorem~\ref{thm:DRkgroups2}.	 Since $j_\Omega$ is injective,  we have
\[
\ker(\tilde\sigma^{-1}\circ \tilde{\imath}\circ\tilde\sigma_H)=\ker(j_\Omega\circ \tilde\sigma^{-1}\circ \tilde{\imath}\circ\tilde\sigma_H)=
\ker(i_*\circ j_H).
\]
Since $j_H$ is injective, we now have
\begin{equation}\label{eq:towardsS}
\ker(i_*)\cap j_H\big(C_c(H,\N) + \image(1-(T_1|_H)_*,1-(T_2|_H)_*)\big) = \{0\}.
\end{equation}
We have $\kappa_H=j_H\circ q$, and
\begin{align*}
   \image(j_H)\cap K_0(C^*(G|_H))_+ &= \kappa_H(C_c(H,\N)) \cap K_0(C^*(G|_H))_+\quad\quad\text{(since $q$ is surjective)}\\
   &=\kappa_H(C_c(H,\N))\quad\quad\text{(using assumption  \eqref{P})}\\
   &= j_H\big(C_c(H,\N) + \image(1-(T_1|_H)_*,1-(T_2|_H)_*)\big).
\end{align*}
	Now \eqref{eq:towardsS} gives
	\[
\ker(i_*)\cap	\image(j_H)\cap K_0(C^*(G|_H))_+ =\{0\}.
	\]
Finally,  consider the right  square of the commuting diagram in Theorem~\ref{thm:DRkgroups2}.  Since $\sigma|^{-1}\circ\imath_*\circ\sigma_H|$ is injective, we have
\[
\image(j_H)=\ker(\tau_H)=\ker (\sigma|^{-1}\circ\imath_*\circ\sigma_H|\circ\tau_H)=\ker(\tau_\Omega\circ i_*)\supset\ker(i_*),
\]	
and hence~\eqref{S} holds.
\end{proof}

	We can now state a simple result about stable finiteness of the $C^*$-algebra of a Deaconu--Renault groupoid with just one nontrivial open invariant subset, involving only condition~\eqref{P} and conditions on the groupoid.

	\begin{corollary}
		Let $T_1$ and $T_2$ be commuting surjective local homeomorphisms on a second-countable totally disconnected locally compact Hausdorff space $\Omega$, and let $G_T$ be the associated rank-$2$ Deaconu--Renault groupoid. Suppose that $\emptyset \subsetneq H \subsetneq \Omega$ is the unique non-trivial open invariant subset of the unit space $\Omega$ of $G_T$.  Suppose that the $T_i$ and the $T_i|_{\Omega \setminus H}$ satisfy~\eqref{M}. If~\eqref{P} holds with respect to $H$, then $C^*(G_T)$ is stably finite.
	\end{corollary}
	\begin{proof} Write $G\coloneqq G_T$.
	To apply Theorem~\ref{thm:lastmain}, we need to verify that $C^*(G|_H)$ and $C^*(G|_{\Omega \setminus H})$ are stably finite. Since $H$ is invariant, it follows from Lemma~\ref{lem:invarianceequivalence} that the groupoid $G|_H$ is equal to the Deaconu--Renault groupoid of the restricted dynamics $T_i|_{H}$.  Lemma~\ref{lem:MforH} implies that the $T_i|_H$ satisfy~\eqref{M}. So, by Proposition~\ref{prop: matrixconditioncoboundary}, the groupoids $G|_H$ and $G|_{\Omega \setminus H}$  are minimal and satisfy the coboundary condition~\eqref{C}. Since $H$ is the unique nontrivial open invariant subset of $\Omega$, these groupoids are both minimal. Hence \cite[Corollary~6.6]{RainoneSims2020} implies that $C^*(G|_H)$ and $C^*(G|_{\Omega \setminus H})$ are stably finite as needed. Now Theorem~\ref{thm:lastmain} gives the result.
	\end{proof}

\appendix\section{}\label{sec:appendix}

This appendix includes the proofs of the technical lemmas, Lemmas~\ref{lem:fcbimodulemorphism}-\ref{lem:fletchermorphisms} used in \S\ref{sec:ktheoryfletcher}, and proofs that certain diagrams commute  used in the proof of Proposition~\ref{prop:square3}.

\subsection{Proofs of Lemmas~\ref{lem:fcbimodulemorphism}--\ref{lem:fletchermorphisms}} \label{ssec:appendix1}

For the proof of Lemma~\ref{lem:fcbimodulemorphism} we use frames, and we now briefly recall what we need to know about them.

Let $X$ be a right-Hilbert $A$-module. A \emph{frame} for $X$ is a sequence $(\eta_i)^\infty_{i=1}$ in $X$ such that $\sum^\infty_{i=1} \theta_{\eta_i, \eta_i}(\xi) = \xi$ for all $\xi \in X$. Equivalently, the partial sums $\sum^n_{i=1} \theta_{\eta_i, \eta_i}$ form an increasing approximate identity for $\K(X)$. If $X$ is countably generated over $A$, then $X$ admits a frame: this can be deduced from the final paragraph of \cite[\S3]{FrankLarson} by regarding $X$ as a countably generated Hilbert module over the unitisation $\widetilde{A}$ of $A$.

Let  $(\eta_m)_m$ be a frame for $X$. Then for all $a\in A$, we can write
	\[\ell(k(a)) = i_1(a) - \sum_m i_1(a)(i_1,i_2)^{(1)}(\theta_{\eta_m,\eta_m}) = i_1(a) - \sum_m i_1(a)i_2(\eta_m)i_2(\eta_m)^*.\]

Lemma~\ref{lem:psimorphismdetails}  is used in the proof of Lemma~\ref{lem:fcbimodulemorphism}.

\begin{lemma}\label{lem:psimorphismdetails}
	Let $\mathbf{Z}$ be a rank-$2$ compactly aligned product system over a $C^*$-algebra $A$, let $\mathcal{N}\TT_{\mathbf{Z}}$ be its Nica--Toeplitz algebra and for $(m,n)\in\N^2$, let $i_{(n,m)}\colon  \mathbf{Z}_{(n,m)}\to \NT_{\mathbf{Z}}$ be the canonical inclusion map. Let $\xi, \zeta \in \mathbf{Z}_{(1,0)}, \eta\in \mathbf{Z}_{(0,1)}$ and $a\in A$. Let $(\eta_m)_m$ be a frame for $\mathbf{Z}_{(0,1)}$. Then
\[
    i_{(0,1)}^{(1)}(\theta_{\eta,\eta})i_{(1,0)}^{(1)}(\theta_{\zeta,\zeta}) i_{(1,0)}(\xi) \Big(i_{(0,0)}(a) - \sum_m i_{(0,0)}(a)i_{(0,1)}(\eta_m)i_{(0,1)}(\eta_m)^*\Big) = 0.
\]
\end{lemma}

\begin{proof}
By Nica covariance,
\[
    i_{(0,1)}^{(1)}(\theta_{\eta,\eta})i_{(1,0)}^{(1)}(\theta_{\zeta,\zeta}) = i_{(1,1)}^{(1)}(T),
\]
for some $T\in \mathcal{K}(\mathbf{Z}_{(1,1)})$. Furthermore, $T$ is in the closed span of elements of the form $\theta_{\alpha\otimes\beta,\gamma\otimes\delta}$ where $\alpha,\gamma \in \mathbf{Z}_{(1,0)}$ and $\beta,\delta\in \mathbf{Z}_{(0,1)}$. Thus it suffices to show that for such $\alpha,\beta,\gamma,\delta$,
\[
    i_{(1,0)}(\alpha)i_{(0,1)}(\beta)i_{(0,1)}(\delta)^*i_{(1,0)}(\gamma)^*i_{(1,0)}(\xi)\Big(i_{(0,0)}(a) - \sum_m i_{(0,0)}(a)i_{(0,1)}(\eta_m)i_{(0,1)}(\eta_m)^*\Big) = 0.
\]
Note that
\[
    i_{(1,0)}(\gamma)^*i_{(1,0)}(\xi) i_{(0,0)}(a) = i_{(0,0)}(\langle \gamma,\xi\rangle)i_{(0,0)}(a) = i_{(0,0)}(\langle \gamma,\xi\rangle a),
\]
and hence it suffices to show that if $b\in A$, then
\[
    i_{(0,1)}(\delta)^*\Big(i_{(0,0)}(b) - \sum_m i_{(0,0)}(b)i_{(0,1)}(\eta_m)i_{(0,1)}(\eta_m)^*\Big) = 0.
\]

There is a copy of $\TT_{\mathbf{Z}_{(0,1)}}$ within $\NT_{\mathbf{Z}}$ and the left hand side of the above equation sits inside this copy: specifically, for the canonical maps $\ell$ and $k$ of \S\ref{ssec:HBktheory},
\[
    i_{(0,1)}(\delta)^*\Big(i_{(0,0)}(b) - \sum_m i_{(0,0)}(b)i_{(0,1)}(\eta_m)i_{(0,1)}(\eta_m)^*\Big)
        = \big(\ell\circ k(b^*)i_2(\delta)\big)^* \in \TT_{\mathbf{Z}_{(0,1)}}.
\]
So it suffices to show that $\ell\circ k(b^*)i_2(\delta) = 0$. Identifying $\TT_{\mathbf{Z}_{(0,1)}}$ with its image in $\LL(\F_{\mathbf{Z}_{(0,1)}})$ via \cite[Proposition~3.3]{Pimsner} as discussed immediately after Notation~\ref{notation:i_1andi_2}, we see that $\ell \circ k(b^*)$ acts on $\F_{\mathbf{Z}_{(0,1)}}$ as multiplication by $b^*$ in the 0 coordinate and zero in the remaining coordinates, while, $i_2(\delta)$ acts by tensoring on the left by $\delta$ and shifting right. Hence
\begin{align*}
    \ell\circ k(b^*)(i_2(\delta)(\sigma_0,\sigma_1,\sigma_2,\dots))
        &= \ell\circ k(b^*)(0,\delta\cdot\sigma_0,\delta\otimes\sigma_1,\delta\otimes \sigma_2,\dots)\\
        &= \ell(b^* 0, 0, 0,\dots) = 0.\qedhere
\end{align*}
\end{proof}

\begin{proof}[Proof of Lemma~\ref{lem:fcbimodulemorphism}]
\label{proof_lem:fcbimodulemorphism}
We first show that there is an isometric linear map $\psi \colon X_1 \to X_1 \otimes_A \mathcal{T}_{X_2}$ satisfying $\psi(\xi \cdot a) = \xi \otimes \ell(k(a))$ for all $\xi \in X_1$ and $a \in A$.
For this, fix $\xi, \zeta \in X_1$ and $a,b\in A$. We have
\begin{align}
\langle \xi \otimes \ell(k(a)),\zeta \otimes \ell(k(b))\rangle & = \langle \ell(k(a)),\langle \xi,\zeta\rangle \cdot \ell(k(b))\rangle \notag \\
	& = \ell(k(a))^*i_1(\langle \xi, \zeta\rangle)\ell(k(b)) \notag \\
	& \overset{*}{=} \ell(k(a^* \langle \xi, \zeta\rangle b)) \notag \\
	& = \ell\circ k(\langle \xi \cdot a,\zeta\cdot b\rangle),  \label{eq:ellkisom}
\end{align} where the equality $*$ holds because of the way the images of $\ell\circ k$ and $i_1$ act on the Fock space. Because $\ell$ and $k$ are injective $*$-homomorphisms and hence isometries, for $\xi_i \in X_1$ and $a_i \in A$,
\begin{align*}
\Big\|\sum_i \xi_i\cdot a_i\Big\|^2
    & = \Big\|\Big\langle \sum_i\xi_i\cdot a_i,\sum_i \xi_i\cdot a_i\Big\rangle \Big\| \\
	& = \Big\| \sum_{i,j} \langle \xi_i\cdot a_i,\xi_j\cdot a_j\rangle \Big\| \\
	& = \Big\|\ell\circ k\Big(\sum_{i,j} \langle \xi_i\cdot a_i, \xi_j \cdot a_j\rangle\Big) \Big\| \\
	& = \Big\|\sum_{i,j} \ell \circ k (\langle \xi_i\cdot a_i, \xi_j \cdot a_j\rangle) \Big\| \\
	& = \Big\|\sum_{i,j} \langle \xi_i\otimes \ell(k(a_i)),\xi_j \otimes \ell(k(a_j))\rangle \Big\| \text{ by \eqref{eq:ellkisom},} \\
	& = \Big\|\Big\langle \sum_i \xi_i \otimes \ell(k(a_i)),\sum_j \xi_j \otimes \ell(k(a_j))\Big\rangle \Big\| \\
    & = \Big\|\sum_i \xi_i\otimes \ell(k(a_i))\Big\|^2.	
\end{align*}
This implies first that if $\sum_i \xi_i\cdot a_i = \sum_j \eta_j\cdot b_j$, then $\|\sum_i \xi_i\otimes \ell(k(a_i)) - \sum_j \eta_i\otimes \ell(k(b_i))\| = 0$, so there is a linear map $\psi \colon \operatorname{span} \{\xi \cdot a : \xi \in X_1, a \in A\} \to X_1 \otimes_A \mathcal{T}_{X_2}$ satisfying $\psi(\xi \cdot a) = \xi \otimes \ell(k(a))$ for all $\xi \in X_1$ and $a \in A$; and second that this map is norm-decreasing. Cohen--Hewitt factorisation \cite[Theorem~2.5]{Hewitt64} gives $X_1 = \{\xi \cdot a : \xi \in X_1, a \in A\}$, so the domain of $\psi$ is all of $X_1$.

Next, we show that $(\ell\circ k, \psi)$ is a covariant right-Hilbert bimodule morphism. Let $a,b\in A$ and $\xi \in X_1$. Then
	\[
        \psi(\xi\cdot a)\cdot \ell\circ k(b) = (\xi\otimes \ell(k(a)))\cdot\ell(k(b)) = \xi\otimes(\ell(k(ab)) = \psi(\xi\cdot(ab)) =\psi((\xi\cdot a)\cdot b),
    \]
proving that $(\ell\circ k, \psi)$ is a right module morphism. That $(\ell\circ k, \psi)$ preserves the inner product, follows from \eqref{eq:ellkisom} by linearity and continuity.
	
To check that the left action is preserved, consider $\ell\circ k(b)\cdot \psi(\xi\cdot a) \in X_1\otimes_A \TT_{X_2}$. Let $\mathbf{Z}$ be the rank-$2$ product system such that $\mathbf{Z}_{(1,0)} = X_1$ and $\mathbf{Z}_{(0,1)} = X_2$.  Since $X_1$ and $X_2$ are countably generated by assumption, there exist frames $(\xi_n)_n$ and $(\eta_m)_m$  for $X_1$ and $X_2$, respectively.
	There is an injection $\gamma : \TT_{X_2} \to \NT_{\mathbf{Z}}$ that carries $\ell\circ k(b)\cdot \psi(\xi\cdot a) \in \TT_{X_2}$ to 	
    \[
        \Xi = \Big(i_{(0,0)}(b) - \sum_m i_{(0,0)}(b)i_{(0,1)}^{(1)}(\theta_{\eta_m,\eta_m})\Big) i_{(1,0)}(\xi) \Big(i_{(0,0)}(a) - \sum_m i_{(0,0)}(a)i_{(0,1)}^{(1)}(\theta_{\eta_m,\eta_m})\Big) \in \NT_{\mathbf{Z}}.
    \]
	Then $i_{(1,0)}(\xi) = \sum_n i_{(1,0)}^{(1)}(\theta_{\xi_n,\xi_n})i_{(1,0)}(\xi)$ and hence by Lemma~\ref{lem:psimorphismdetails},
\begin{align*}
    \Xi & = i_{(0,0)}(b) i_{(1,0)}(\xi) \Big(i_{(0,0)}(a) - \sum_m i_{(0,0)}(a)i_{(0,1)}^{(1)}(\theta_{\eta_m,\eta_m})\Big)\\
		& = i_{(1,0)}(b\cdot\xi) \Big(i_{(0,0)}(a) - \sum_m i_{(0,0)}(a)i_{(0,1)}^{(1)}(\theta_{\eta_m,\eta_m})\Big),
\end{align*}
which is $\gamma(i_1(b\cdot \xi) \ell\circ k(a)) \in \gamma(\TT_{X_2}) \subset \NT_{\mathbf{Z}}$. Thus $\ell\circ k(b)\cdot \psi(\xi\cdot a) =\psi((b\cdot\xi)\cdot a)=\psi(b\cdot(\xi\cdot a))$, as needed.

Finally for covariance, note that $A$ acts compactly and injectively on the right-Hilbert $\mathcal{T}_{X_2}$-module $\mathcal{T}_{X_2}$, and hence the map $l \circ k$ carries $\phi_{X_1}^{-1}\big(\mathcal{K}(X_1) \cap (\ker \phi_{X_1})^\perp\big)$ into $(\phi_{X_1} \otimes 1_{\TT_{X_2}})^{-1} \big(\mathcal{K}(X_1 \otimes_A \TT_{X_2}) \cap (\ker (\phi_{X_1} \otimes 1_{\TT_{X_2}}))^\perp\big)$.
\end{proof}

\begin{proof}[Proof of Lemma \ref{lem:incbimodulemorphism}]
\label{proof_lem:incbimodulemorphism}
We proceed as in the proof  of Lemma~\ref{lem:fcbimodulemorphism}. For $\xi,\zeta \in X_1$ and $a,b \in A$ we have
	\begin{align*}
		\langle \xi \otimes i_1(a),\zeta \otimes i_1(b)\rangle & = \langle i_1(a),\langle \xi, \zeta\rangle \cdot i_1(b) \rangle \\
		& = i_1(a)^* i_1(\langle \xi, \zeta \rangle)i_1(b) \\
		& = i_1(a^*\langle \xi, \zeta\rangle b) = i_1(\langle \xi\cdot a, \eta \cdot b\rangle).
	\end{align*}
Since $i_1\colon  A\to \TT_{X_2}$ is an isometry, it follows as in the proof of  Lemma~\ref{lem:fcbimodulemorphism} that there is an isometric linear map $\varphi \colon X_1 \to X_1 \otimes_A \mathcal{T}_{X_2}$ such that $\varphi(\xi \cdot a) = \xi \otimes i_1(a)$ for all $\xi \in X_1$ and $a \in A$. That $\varphi$ is a right module morphism and intertwines inner-products is trivial. Since the left action of $\TT_{X_2}$ on $X_1 \otimes \TT_{X_2}$, restricted to $A$, is just multiplication by $i_1(A)$, we see that $\varphi$ is a left module morphism. The proof of covariance is as in Lemma~\ref{lem:fcbimodulemorphism}.
\end{proof}

\begin{proof}[Proof of Lemma \ref{lem:fletcherisomorphism}]
\label{proof_lem:fletcherisomorphism}
    We use Fletcher's commuting diagram~\cite[Figure~3.9]{FletcherThesis} in the instance where $P = Q = \mathbb{N}$ and $\mathbf{Z}_{p,q} = X_1^{\otimes p} \otimes_A X_2^{\otimes q}$ for $p,q \in \mathbb{N}$.
    There is an isomorphism of $X_1 \otimes_A \mathcal{T}_{X_2}$ with Fletcher's $Y^{\mathcal{NT}}_{(1, 0)}$ that carries $\xi \otimes i_1(a)$ to $i_{\mathbf{Z}_{(1,0)}}(\xi)\phi^{\mathcal{NT}_\mathbf{X}}(a)$ in Fletcher's notation, and there is a similar isomorphism
    of $\mathbf{W}^{\mathcal{NO}}_{(0,1)}$ with $X_2 \otimes_A \mathcal{O}_{X_1}$. Making these replacements in Fletcher's diagram, the entries $\mathcal{NO}_{\mathbf{Y}^\mathcal{NT}}$ and $\mathcal{NT}_{\mathbf{W}^\mathcal{NO}}$
    in his diagram become  $\mathcal{O}_{X_1 \otimes \mathcal{T}_{X_2}}$ and $\mathcal{T}_{X_2 \otimes \mathcal{O}_{X_1}}$ since the Nica--Toeplitz and Cuntz--Nica--Pimsner algebra of a rank-1 product system are just the standard Toeplitz algebra and Cuntz--Pimsner algebra of its
    generating fibre \cite[Proposition~2.11]{Fowler}. With these identifications in place, let $J$ be the map $\omega$ from top right to top left in Fletcher's diagram. Fix $\xi \in X_1$ and $a \in A$. Then $j_1(\xi \otimes i_1(a))$ is the image of
    $\xi \cdot a \in \mathbf{Z}_{(1,0)}$ under the composite map $\xi \cdot a \mapsto \xi \otimes_A i_1(a) \mapsto j_1(\xi \otimes i_1(a)))$, which is the composite map from $\mathbf{Z}_{(1,0)}$ to $\mathcal{O}_{X_1 \otimes_A \mathcal{T}_{X_2}}$ in Fletcher's diagram.
    Likewise, $i_1(j_2(\xi \cdot a))$ is the image of the same element under the composite map from $\mathbf{Z}_{(1,0)}$ to $\mathcal{T}_{X_2 \otimes \mathcal{O}_{X_1}}$ in Fletcher's diagram. So $J(j_1(\xi \otimes i_1(a)) = i_1(j_2(\xi \cdot a))$ because Fletcher's diagram commutes.

  Now  take $(p,q)=(0,0)$ so that $\mathbf{Z}_{0,0} =A$.   Then  the homomorphisms $i_{\mathbf{Z}_{0,0}}$ and $j_{Y_0}^\mathcal{NT}$ of Fletcher's diagram are our $j_1\colon  \mathcal{T}_{X_2}\to \mathcal{O}_{X_1\otimes\mathcal{T}_{X_2}}$ and  $i_1\colon  A\to \mathcal{T}_{X_2}$, respectively, and  the homomorphisms $j_{\mathbf{Z}_{(0,0)}}$ and $i_{\mathbf{W}_0^{\mathcal{NO}}}$ are our $j_1\colon A\to \mathcal{O}_{X_1}$ and $i_1\colon  \mathcal{O}_{X_1}\to \mathcal{T}_{X_2\otimes \mathcal{O}_{X_1}}$, respectively.  Thus it follows from Fletcher's diagram that  $J(j_1(i_1(a))) = i_1(j_1(a))$.

 Similarly, for $\eta' \in X_2$, the element $j_1(i_2(\eta' \cdot \langle \eta', \eta'\rangle))$ is the image of $\eta' \cdot \langle \eta', \eta'\rangle$ under the composite map from $\mathbf{Z}_{(0,1)}$ to $\mathcal{O}_{X_1 \otimes_A \mathcal{T}_{X_2}}$  in Fletcher's diagram, while $i_2(\eta' \otimes j_1(\langle \eta', \eta'\rangle))$ is the image of the same element under the composite map from $\mathbf{Z}_{(1,0)}$ to $\mathcal{T}_{X_2 \otimes \mathcal{O}_{X_1}}$ in Fletcher's diagram. So once again commutativity of Fletcher's
    diagram implies that $J(j_1(i_2(\eta' \cdot \langle \eta', \eta'\rangle))) = i_2(\eta' \otimes j_1(\langle \eta', \eta'\rangle))$.
   \end{proof}

\begin{proof}[Proof of Lemma \ref{lem:fletchermorphisms}]
\label{proof_lem:fletchermorphisms}
	The first statement is trivial. In the second statement, the map $1\otimes \kappa_1$ is well defined because $\kappa_1$ is induced by the bimodule morphism of the first statement. That the second pair is a right module morphism is obvious. Given $\xi,\eta \in X_2I$ and $a,b \in \OO_{X_1I}$,
\begin{align*}
    \langle 1\otimes\kappa_1(\xi\otimes a),1\otimes \kappa_1(\eta\otimes b)\rangle
        &= \langle \xi\otimes\kappa_1(a),\eta\otimes\kappa_1(b)\rangle \\
        &= \langle \kappa_1(a),\langle \xi,\eta\rangle \cdot \kappa_1(b)\rangle \\
		&= \kappa_1(a)^*j_1(\langle \xi,\eta\rangle)\kappa_1(b) \\
        &= \kappa_1(a)^*\kappa_1(j_1(\langle \xi,\eta\rangle))\kappa_1(b) \\
        &= \kappa_1(a^*j_1(\langle \xi,\eta\rangle)b) \\
        &= \kappa_1(\langle \xi\otimes a,\eta \otimes b\rangle).
\end{align*}
Thus, by linearity and continuity, the inner product condition for a Hilbert-module morphism is satisfied.
    To see that $(\kappa_1, 1 \otimes \kappa_1)$ is a left-module morphism, since the left action of the coefficient algebra on a Hilbert module is implemented by a $C^*$-homomorphism, it suffices to show that the left actions of generators of $\OO_{X_1 I}$ on $Y_1^I$ are preserved by $\kappa_1$ and $1 \otimes \kappa_1$. That is, for $a \in I$, $x_1 \in X_1I$, $x_2 \in X_2 I$ and $b \in \OO_{X_1 I}$,  we must show that
    \begin{align*}
    \kappa_1(j_1(a))\big((1 \otimes \kappa_1)(x_2 \otimes b)) &= (1 \otimes \kappa_1)(j_1(a) \cdot (x_2 \otimes b))
        \qquad\text{ and }\\
    \kappa_1(j_2(x_1))\big((1 \otimes \kappa_1)(x_2 \otimes b)) &= (1 \otimes \kappa_1)(j_2(x_1) \cdot (x_2 \otimes b)).
    \end{align*}
    For the first of these equations, we just calculate
    \[
    \kappa_1(j_1(a))\big((1 \otimes \kappa_1)(x_2 \otimes b)) = j_1(a)\cdot (x_2 \otimes b) = (a \cdot x_2) \otimes b = (1 \otimes \kappa_1)(j_1(a) \cdot (x_2 \otimes b)).
    \]

    For the second equation, recall that $X_1 \otimes_{A} X_2 \cong \mathbf{X}_{(1,1)} \cong X_2 \otimes_{A} X_1$, and this isomorphism extends the corresponding isomorphism $X_1 I \otimes_{I} X_2 I \cong \mathbf{XI}_{(1,1)} \cong X_2 I \otimes_{I} X_1 I$.
    So there is a sequence, indexed by $\ell \in \N$, of finite linear combinations $\sum_{k=1}^{K_\ell} x_2^{k, \ell} x_1^{k, \ell}$ in $\mathbf{XI}_{(1,1)}$ with each $x^{k, \ell}_i \in X_i I$ such that
    $\sum_{k=1}^{K_\ell} x_2^{k, \ell} x_1^{k, \ell} \to x_1x_2$. By definition (see \cite[Lemma~3.1.22 and Proposition~3.3.6]{FletcherThesis}),
    \[
        j_2(x_1)\cdot(x_2 \otimes b) = \lim_\ell \sum^{K_\ell}_{k=1} x_2^{k, \ell} \otimes j_2(x_1^{k, \ell})b.
    \]
    Since exactly the same approximation holds in $X_1 \otimes_{A} X_2 \cong \mathbf{X}_{(1,1)} \cong X_2 \otimes_{A} X_1$, the left action of $\kappa_1(j_2(x_1))$ on $(1 \otimes \kappa_1)(x_2 \otimes b)$ is given by the same formula. Again, the proof of covariance is as in Lemma~\ref{lem:fcbimodulemorphism}.
\end{proof}

\subsection{Results used in the proof of Proposition~\ref{prop:square3}}

\begin{lemma}\label{lem:subgdiagrama}
	Let $\eta \colon  \K(\F_{X_1})\to \K(\F_{X_1\otimes \TT_{X_2}})$ be the homomorphism induced by the right-Hilbert bimodule morphism $\psi \colon  X_1 \to X_1 \otimes_A \TT_{X_2}$ of Lemma~\ref{lem:fcbimodulemorphism}, so that $\psi(\xi \cdot a)=\xi \otimes \ell(k(a))$ for $\xi \in X_1, a \in A$.
	The diagram
	\[
	\begin{tikzcd}
		A \arrow[d, "k"] \arrow[r, "k"]         & {\K(\F_{X_2})} \arrow[r, "\ell"] & {\TT_{X_2}} \arrow[d, "k"]            \\
		{\K(\F_{X_1})} \arrow[rr, "\eta"] &                                        & {\K(\F_{X_1\otimes \TT_{X_2}}),}
	\end{tikzcd}\] from the proof of Proposition \ref{prop:square3} commutes.
\end{lemma}

\begin{proof}
	First, we describe what $\eta$ is. For each $n\ge 1$, $\psi$ induces a map $\psi_n \colon  X_1^{\otimes n} \to (X_1 \otimes_A \TT_{X_2})^{\otimes n}$ such that $\psi_n(\xi_1\otimes\cdots\otimes\xi_n) = \psi(\xi_1)\otimes\cdots\otimes \psi(\xi_n)$. These $\psi_n$ together with the map $\psi_0 := \ell\circ k$ induce a homomorphism $\Psi = \bigoplus_{n=0}^\infty \psi_n \colon  \mathcal{F}_{X_1} \to \mathcal{F}_{X_1\otimes \TT_{X_2}}$ which in turn induces $\eta \coloneqq \Psi^{(1)} \colon  \mathcal{K}(\mathcal{F}_{X_1})\to \mathcal{K}(\mathcal{F}_{X_1\otimes \TT_{X_2}})$. Now let $a,b \in A$ and $\zeta \in X^{\otimes n}$ (considered as an element of the Fock space) for $n \ge 0$. Starting with $ab^* \in A$ and tracing the diagram along the top and right, we get
	\[k(\ell(k(ab^*)))(\zeta) = \begin{cases}
		\big(\ell(k(ab^*))\zeta,0,0,\dots\big), & \text{ if } n=0, \\
		0, \text{ if } n\ge 1.
	\end{cases}\]
	Tracing the diagram along the left and bottom, we instead get
	\begin{align*}
		\eta(k(ab^*))(\zeta) & = \Psi^{(1)}(\theta_{(a,0,0,\dots),(b,0,0,\dots)})(\zeta) \\
		& = \theta_{\Psi(a,0,0,\dots),\Psi(b,0,0,\dots)}(\zeta) \\
		& = \theta_{\left(\psi_0(a),0,0,\dots\right),\left(\psi_0(b),0,0,\dots\right)}(\zeta) \\
		& = \left(\psi_0(a),0,0,\dots\right)\cdot \langle \left(\psi_0(b),0,0,\dots\right),\zeta\rangle \\
		& = \begin{cases}
			\big(\ell(k(a)),0,0,\dots\big)\cdot \langle \ell(k(b)), \zeta\rangle, & \text{ if } n=0,\\
			0, & \text{ if } n\ge 1,
		\end{cases} \\
		& = \begin{cases}
			\big(\ell(k(a)),0,0,\dots\big)\cdot(\ell(k(b))^*\zeta), & \text{ if } n=0,\\
			0, & \text{ if } n\ge 1,
		\end{cases} \\
	& = \begin{cases}
		\big(\ell(k(ab^*))\zeta ,0,0,\dots\big), & \text{ if } n=0,\\
		0, & \text{ if } n\ge 1.
	\end{cases}
	\end{align*} Thus the diagram commutes.
\end{proof}

\begin{lemma}\label{lem:subgdiagramb}
	The diagram
	\[
	\begin{tikzcd}
		{\OO_{X_1}} \arrow[d, "\theta"] \arrow[r, "k"]   & {\K(\F_{X_2\otimes\OO_{X_1}})} \arrow[d, "\ell"] \\
		{\OO_{X_1\otimes\TT_{X_2}}} \arrow[r, "J"] & {\TT_{X_2\otimes\OO_{X_1}},}
	\end{tikzcd}\] from the proof of Proposition \ref{prop:square3} commutes.
\end{lemma}
\begin{proof}
Fix $a\in A$ and $\xi \in X_1$ so that $j_2(\xi\cdot a) \in \OO_{X_1}$. Tracing the diagram along the top and right gives us $\ell(k(j_2(\xi\cdot a))) = (j_2(\xi\cdot a),0,0,\dots)$ as an operator on the Fock space. Since $X_2$ is countably generated, we can choose a frame  $(\eta_m)_m$  for $X_2$. If instead we go down the left side, first note that
\begin{align*}
    \theta(j_2(\xi\cdot a)) = j_2(\xi\otimes \ell(k(a)))
        & = j_2(\xi \otimes i_1(a)) - \sum_m j_2(\xi \otimes i_1(a)i_2(\eta_m)i_2(\eta_m)^*) \\
        & = j_2(\xi \otimes i_1(a)) - \sum_m j_2\left(\left(\xi \otimes i_1(a) \right)\cdot(i_2(\eta_m)i_2(\eta_m)^*)\right) \\
        & = j_2(\xi \otimes i_1(a)) - \sum_m j_2\left(\xi \otimes i_1(a) \right)j_1(i_2(\eta_m))j_1(i_2(\eta_m))^*.
\end{align*}
	For each $m$, write $\eta_m = \eta_m'\cdot \langle \eta_m',\eta_m'\rangle$. Then by Lemma \ref{lem:fletcherisomorphism},
	\begin{align*}
		J(\theta(j_2(\xi \cdot a))) = i_1(j_2(\xi\cdot a)) - i_1(j_2(\xi\cdot a))i_2(\eta_m'\otimes j_1(\langle \eta_m',\eta_m'\rangle))i_2(\eta_m'\otimes j_1(\langle \eta_m',\eta_m'\rangle))^*.
	\end{align*}
	Now if $a_0 \in \OO_{X_1} = (X_2\otimes\OO_{X_1})^{\otimes 0},$ then \[i_2(\eta_m'\otimes j_1(\langle \eta_m',\eta_m'\rangle))i_2(\eta_m'\otimes j_1(\langle \eta_m',\eta_m'\rangle))^*(a_0) = 0,\] since for any Toeplitz algebra $\TT_Y$ and $y\in Y$, $i_2(y)^*$ sends a vector $(c_0,0,0,\dots)$ in the Fock space to $0$. Furthermore, given any $\eta \in X_2$ and $b\in \OO_{X_1}$ so that $\eta \otimes b \in X_2\otimes_A \OO_{X_1}$, we have by \cite[Lemma 1.9]{Katsura} that if $m\ge 1$ and $\bigotimes_{n=1}^{m} (\zeta_n\otimes a_n) \in (X_2\otimes_A\OO_{X_2})^{\otimes m}$,
\begin{align*}
	i_2(\eta\otimes b)i_2(\eta\otimes b)^*\Big(\bigotimes_{n=1}^{m} (\zeta_n\otimes a_n)\Big)
        & = \theta_{\eta\otimes b,\eta\otimes b}(\zeta_1\otimes a_1)\otimes \bigotimes_{n=2}^{m}(\zeta_n\otimes a_n) \\
		& = \left(\eta\otimes b \cdot \langle \eta \otimes b, \zeta_1\otimes a_1\rangle\right)\otimes \bigotimes_{n=2}^{m}(\zeta_n\otimes a_n) \\
		& = \left(\eta \otimes b \cdot \left (b^* j_1(\langle \eta,\zeta_1\rangle)a_1 \right)\right)\otimes \bigotimes_{n=2}^{m}(\zeta_n\otimes a_n) \\
		& = \left(\eta \otimes \left(bb^*j_1(\langle \eta,\zeta_1\rangle)a_1\right)\right)\otimes \bigotimes_{n=2}^{m}(\zeta_n\otimes a_n).
\end{align*}
Taking $\eta = \eta_m'$ and $b = j_1(\langle \eta_m',\eta_m'\rangle)$, we have that
\begin{align*}
	i_2(\eta_m'\otimes j_1(\langle \eta_m',\eta_m'\rangle))&i_2(\eta_m'\otimes j_1(\langle \eta_m',\eta_m'\rangle))^*\Big(\bigotimes_{n=1}^{m} (\zeta_n\otimes a_n)\Big) \\
		& = \left(\eta_m' \otimes (j_1(\langle \eta_m',\eta_m'\rangle)j_1(\langle \eta_m',\eta_m'\rangle)^* j_1(\langle \eta_m',\zeta_1\rangle)a_1)\right)\otimes \bigotimes_{n=2}^{m}(\zeta_n\otimes a_n) \\
		& = \left(\left(\eta_m' \cdot \langle \eta_m',\eta_m'\right) \otimes \big(j_1(\langle\eta_m'\cdot\langle\eta_m',\eta_m'\rangle, \zeta_1\rangle) a_1\big) \right)\otimes \bigotimes_{n=2}^{m}(\zeta_n\otimes a_n) \\
		& = \left(\eta_m\otimes \big(j_1(\langle \eta_m,\zeta_1\rangle) a_1\big)\right)\otimes \bigotimes_{n=2}^{m}(\zeta_n\otimes a_n) \\
		& = \left(\left(\eta_m\cdot\langle \eta_m,\zeta_1\rangle\right)\otimes a_1\right)\otimes \bigotimes_{n=2}^{m}(\zeta_n\otimes a_n).
\end{align*}
Thus
\begin{align*}
	\sum_m i_2(\eta_m'\otimes j_1(\langle \eta_m',\eta_m'\rangle))
        & i_2(\eta_m'\otimes j_1(\langle \eta_m',\eta_m'\rangle))^*\Big(\bigotimes_{n=1}^{m} (\zeta_n\otimes a_n)\Big) \\
        & = \Big(\Big(\sum_m \eta_m \cdot \langle \eta_m,\zeta_1\rangle\Big)\otimes a_1\Big)\otimes \bigotimes_{n=2}^{m}(\zeta_n\otimes a_n) \\
		& = (\zeta_1\otimes a_1) \otimes \bigotimes_{n=2}^{m}(\zeta_n\otimes a_n)
          = \bigotimes_{n=1}^{m}(\zeta_n\otimes a_n).
	\end{align*}
	
	In particular, for any $\zeta = \bigotimes_{n=1}^m(\zeta_n\otimes a_n) \in (X_2 \otimes \OO_{X_1})^{\otimes m}$,
	\[
    J(\theta(j_2(\xi \cdot a)))\Big(\bigotimes_{n=1}^{m}(\zeta_n\otimes a_n)\Big)
        = \begin{cases}
            j_2(\xi\cdot a)\zeta - 0 = j_2(\xi\cdot a)\zeta, & \text{ if } m = 0 \text{ (and hence $\zeta\in \OO_{X_1}$) }, \\
        	j_2(\xi\cdot a)\zeta - j_2(\xi\cdot a)\zeta = 0, & \text{ if } m\ge 1,
		  \end{cases}
    \]
    and thus $J(\theta(j_2(\xi \cdot a))) = (j_2(\xi\cdot a),0,0,\dots)$ and the diagram commutes when beginning with elements of the form $j_2(\xi\cdot a)$. Proving that the diagram commutes when beginning with elements of the form $j_1(a)$ follows a (simpler) similar argument.
\end{proof}

\begin{lemma}\label{lem:subgdiagramc}
	The diagram
	\[
	\begin{tikzcd}
		A \arrow[d, "i_1"] \arrow[r, "k"] & {\K(\F_{X_1})} \arrow[d, "\rho"]     \\
		{\TT_{X_2}} \arrow[r, "k"]  & {\K(\F_{X_1\otimes\TT_{X_2}}),}
	\end{tikzcd}\] from the proof of Proposition \ref{prop:square3} commutes.
\end{lemma}

\begin{proof}
	As in the proof of Lemma \ref{lem:subgdiagrama}, we describe what $\rho$ is. For each $n\ge 1$, $\varphi$ induces a map $\varphi_n \colon  X_1^{\otimes n} \to (X_1 \otimes_A \TT_{X_2})^{\otimes n}$ such that $\varphi_n(\xi_1\otimes\cdots\otimes\xi_n) = \varphi(\xi_1)\otimes\cdots\otimes \varphi(\xi_n)$. These $\varphi_n$ together with the map $\varphi_0 \coloneqq i_1$ induces a homomorphism $\Phi = \bigoplus_{n=0}^\infty \varphi_n \colon  \mathcal{F}_{X_1} \to \mathcal{F}_{X_1\otimes \TT_{X_1}}$ which in turn induces $\rho \coloneqq \Phi^{(1)} \colon  \mathcal{K}(\mathcal{F}_{X_1})\to \mathcal{K}(\mathcal{F}_{X_1\otimes \TT_{X_2}})$. Now let $a,b \in A$ and $\zeta \in X^{\otimes n}$ (considered as an element of the Fock space) for $n \ge 0$. Starting with $ab^* \in A$ and tracing the diagram along the left and bottom, we get
	\[k(i_1(ab^*))(\zeta) = \begin{cases}
		\big(i_1(ab^*)\zeta,0,0,\dots\big), & \text{ if } n=0, \\
		0, \text{ if } n\ge 1.
	\end{cases}\]
	Tracing the diagram along the top and right, we instead get
	\begin{align*}
		\rho(k(ab^*))(\zeta) & = \Phi^{(1)}(\theta_{(a,0,0,\dots),(b,0,0,\dots)})(\zeta) \\
		& = \theta_{\Phi(a,0,0,\dots),\Phi(b,0,0,\dots)}(\zeta) \\
		& = \theta_{\left(\varphi_0(a),0,0,\dots\right),\left(\varphi_0(b),0,0,\dots\right)}(\zeta) \\
		& = \left(\varphi_0(a),0,0,\dots\right)\cdot \langle \left(\varphi_0(b),0,0,\dots\right),\zeta\rangle \\
		& = \begin{cases}
			(i_1(a),0,0,\dots)\cdot \langle i_1(b), \zeta\rangle, & \text{ if } n=0,\\
			0, & \text{ if } n\ge 1,
		\end{cases} \\
		& = \begin{cases}
			(i_1(a),0,0,\dots)\cdot(i_1(b)^*\zeta), & \text{ if } n=0,\\
			0, & \text{ if } n\ge 1,
		\end{cases} \\
		& = \begin{cases}
			(i_1(ab^*)\zeta ,0,0,\dots), & \text{ if } n=0,\\
			0, & \text{ if } n\ge 1.
		\end{cases}
	\end{align*} Thus the diagram commutes.
\end{proof}

\end{document}